\documentclass[11pt]{article}
\usepackage{amsfonts}
\usepackage{mathrsfs}
\usepackage{amssymb}
\usepackage{amsmath}
\usepackage{bbm}
\usepackage{cite}
\usepackage{graphicx}
\usepackage{color}
\usepackage{extarrows}

\newtheorem{theorem}{Theorem}[section]
\newtheorem{prop}[theorem]{Proposition}
\newtheorem{cor}[theorem]{Corollary}
\newtheorem{lemma}[theorem]{Lemma}
\newtheorem{define}[theorem]{Definition}

\def\no{\nonumber}

\newcommand\btd{\raise 2pt \hbox{$\hat\bigtriangledown$}\hskip 1.5pt}
\newcommand\bt{\raise 2pt \hbox{$\bigtriangledown$}\hskip 1.5pt}

\def\no{\nonumber}
\def\x{\textbf{x}}
\def\u{\textbf{u}}
\def\F{\mathcal {F}}
\def\H{\mathcal {H}}
\begin{document}
\title{Symmetry structure of multi-dimensional time-fractional partial differential equations}
\author{Zhi-Yong Zhang \footnote{E-mail:
zzy@muc.edu.cn}\ \ \ \ \ \ \ Jia Zheng \footnote{E-mail:
zhengjia2014@muc.edu.cn}% \ \ \ \ \ \
 \\
\small  College of Science, Minzu University of China,
Beijing 100081, P.R. China}
\date{}
\maketitle
\noindent{\bf Abstract:} In this paper, we concentrate on the Lie symmetry structure of a system of multi-dimensional time-fractional partial differential equations (PDEs). Specifically, we first give an explicit prolongation formula involving Riemann-Liouville time-fractional derivative for the Lie infinitesimal generator in multi-dimensional case, and then show that the infinitesimal generator has an elegant structure. Furthermore, we present two simple conditions to determine the  infinitesimal generators where one is a system of linear time-fractional PDEs, the other is a system of integer-order PDEs and plays the dominant role in finding the infinitesimal generators. We study three time-fractional PDEs to illustrate the {efficiencies} of the results.

\section{Introduction}
The theory of fractional calculus goes back to the Leibniz's letter to L'Hospital \cite{lei}. With the rapid development and extensive applications in the last several decades, nowadays fractional PDEs take an important position in describing the phenomena in the fields such as physics, biology and chemistry, where under certain circumstances  integer-order PDEs cannot work well \cite{ki-1993,pod-1999,sun-2018}. For instance, the anomalous diffusion of tracer particles in complex liquids where superdiffusion and subdiffusion occur is more suitable to be described by the fractional PDE \cite{mc-2014}. {Fractional-order dynamics of the fractional Bloch-Torrey equation, a generalization of the Bloch-Torrey equation by incorporating a fractional order Brownian model of diffusivity, is observed to fit the signal attenuation in diffusion-weighted images obtained from human articular cartilage and human brain \cite{rl-2008}.} Consequently, considerable attentions have been paid to study fractional PDEs and thus a number of effective techniques and methods have been proposed \cite{sri-2006,guo-2015,lu-2010,chen-2008,ahma-2015}.

Lie group theory provides widely applicable techniques to study integer-order PDEs, for example, constructing similarity solutions and {linearized} mappings, investigating integrability, analyzing stability and global {behaviors} of solutions,  etc \cite{olv,lv-1982,blu,he-1994}.
Concerning Lie group theory for fractional PDEs, Buckwar and Luchko first established the invariance of a linear fractional diffusion equation describing subdiffusion in the fractal time random walk \cite{by-1998}
\begin{equation}\label{diff}
\partial{_t^\alpha}u=D\,u_{xx}, ~~~~\alpha>0,
\end{equation}
under the scaling transformations $\bar{x}=\lambda x,\,\bar{t}=\lambda^{2/\alpha}t,\,\bar{u}=u$, where $\partial{_t^\alpha}$ denotes the Riemann-Liouville fractional derivative,  $\lambda$ is a parameter and $D$ is the constant diffusion coefficient. Consequently, by means of the similarity reduction technique \cite{olv}, Eq.(\ref{diff}) was transformed into a fractional ordinary differential equation involving the Erd\'{e}lyi-Kober differential operator and its solutions {were} expressed by the generalized Wright
functions. Gazizov et al. performed Lie symmetry analysis for the nonlinear time-fractional diffusion equations with variable diffusion coefficient $k(u)$
\begin{equation}
\no \partial{_t^\alpha}u=(k(u)u_{x})_x~~~~0<\alpha\leq2,
\end{equation}
{in the sense of} the Riemann-Liouville and Caputo fractional derivatives respectively. The results {showed} that the admitted symmetries in both {fractional-order} cases are narrower than the ones of integer-order case due to the {effects} of time-fractional {derivatives} but still exert important roles in constructing exact solutions and studying symmetry properties \cite{gaz-2007,gaz-2009}.

Following the established schemas symmetry classifications, symmetry reductions and similarity solutions {were} performed for numerous scalar fractional PDEs where such fractional PDEs {originated} from either the descriptions of natural phenomenon \cite{sun-2018,guo-2015} or the direct deformations from the celebrated mathematical physics equations such as the fifth-order KdV equation, Sharma-Tasso-Olver equation, Harry-Dym equation, etc. \cite{huang-2014,liu-stu,wang-2014,zhang-2018,sg-2017,chen-2017}. In \cite{jer-2014}, Jefferson and Carminati wrote an automated package to compute Lie {symmetries} of fractional differential equations under an assumption which was shown to be correct in \cite{zhang-2019}.
In addition to the lower-dimensional scalar time-fractional PDEs, in fact, multi-dimensional fractional PDEs also took effective roles in describing abnormal behaviors \cite{sun-2018}, but only { a small number of papers extended} Lie group theory to study certain special multi-dimensional time-fractional PDEs.
For example, Leo et.al employed Lie group theory for $(1+N)$-dimensional fractional PDEs and studied a scalar fractional diffusion-type equation describing the diffusion of charged particle in a magnetic field \cite{rosa-2017}. Several coupled time-fractional PDEs artificially deformed from classical PDEs were investigated by Lie group theory and affluent exact solutions {were} constructed \cite{ss-2018,pra-2017,kd-2018}. In \cite{zhang-2020}, we used Lie symmetry method for a $(1+2)$-dimensional time-fractional biological population model which describes the changes of population density at the concerned region and found several exact solutions. Such results further demonstrate that Lie symmetry method is a powerful technique to study fractional PDEs.

The prerequisite of applying Lie symmetry method is { the fractional PDEs having} affluent symmetries, but the determining system of Lie symmetries for fractional PDEs contains {the operations of fractional integral and derivative and integer-order derivative}, moreover, the size of the determining system is very large, thus it is not easy to find solutions {of the determining system}. Such dilemma motivates us to explore new techniques to simplify the determining system. Observe that knowing the symmetry structure in advance will greatly drop off the scale of the symmetry determining equations and further facilitate the equation solving \cite{blu,olv}. {Quite recently in \cite{zhang-2019}, by analyzing the structure of the symmetry determining {conditions based on the independence of} time-fractional integrals and derivatives, we {showed} that the infinitesimal generators of Lie symmetries for a scalar time-fractional PDE possess a simple and unified expression. Furthermore, the infinitesimal generators are completely determined by two conditions where one is a system of linear time-fractional PDEs and the other is a system of integer-order PDEs, which makes the Lie symmetries of the scalar fractional PDEs more easier to be found and also the procedure of finding Lie symmetries more convenient to be performed with the known solvers of integer-order PDEs.}
Therefore, as the development of Lie group theory from the lower-dimensional scalar fractional PDEs to the system of multi-dimensional PDEs,  it is significant to systematically and profoundly study Lie group theory of multi-dimensional fractional PDEs and make clear the general symmetry information.
 %Though the infinitesimal generators of some particular multi-dimensional fractional PDEs take the similar form, the general form of the infinitesimal generator is still unknown.there exist no systematic and general results.

In this paper, %we extend the results into multi-dimensional case in the sense of Riemann-Liouville fractional partial derivative. More specifically,
we further investigate Lie symmetry structure of the system consisting of $q$ multi-dimensional time-fractional PDEs with $k$-th order in the sense of Riemann-Liouville fractional derivative, and briefly denoted by
\begin{equation}\label{eqn}
\partial{_t^\alpha}\u=\mathcal{E}(t,\x,\u^{(k)}),~~~0<\alpha<1,
\end{equation}
where $\mathcal {E}=(\mathcal {E}_1, \mathcal {E}_2,\dots, \mathcal{E}_q)$ is a smooth vector function involving $p$ independent variables
$\textbf{x}=(x_1,\dots, x_p)\in \mathbb{R}^p$ and $q$ dependent variables $\textbf{u}=(u_1,\dots, u_q )\in \mathbb{R}^q$,
together with the derivatives of $u_{s}$ with respect to the $x_i\, (s = 1, \dots, q; i = 1, \dots, p)$ up to some order $k$,
denoted by $\textbf{u}^{(k)}=\{u_s^\theta,\,s=1,\dots,q,|\theta|\leq k\}$ with $u_s^\theta=\partial^{|\theta|} u_s/\partial x_1^{\theta_1}\dots\partial x_p^{\theta_p}$, $\theta=\{\theta_1,\dots,\theta_p\}\in Z_{+}^p$ ($Z_{+}$ is the nonnegative integer set) and $|\theta|=\theta_1+\dots+\theta_p\leq k$.
In order to facilitate the analysis of { the symmetry structure}, according to whether the terms in right side of system (\ref{eqn}) are independent of $\u$ and its $\x$-derivatives, we rearrange it as the following form
\begin{equation}\label{eqn-ger}
\partial{_t^\alpha}\u=\F(t,\x,\u^{(k)})+\H(t,\x),~~~0<\alpha<1,
\end{equation}
where $\F=(\F_1, \F_2,\dots, \mathcal{F}_q)$  and $\H=(\H_1, \H_2,\dots, \mathcal
{H}_q)$ are two vector functions, {$\mathcal{E}_i=\mathcal{F}_i+\mathcal{H}_i$, each $\F_i$ collects all the terms containing at least one element of the set $\{u_1,\dots,u_q\}\cup \textbf{u}^{(k)}$ while the remainders in $\mathcal{E}_i$ are collected in $\H_i$ which is only a function of $t$ and $\x$.}

{The main contribution of the paper is to show that the infinitesimal generators of Lie symmetries of multi-dimensional system (\ref{eqn-ger}) have a simple and unified form and are completely determined by two conditions similar as the ones of scalar time-fractional PDE. The two key points of achieving the goals are first to find an explicit prolongation formula  of the infinitesimal generator {in multi-dimensional fractional case} and then to figure out the structure of $\mu_s$ in (\ref{mu}). It should be pointed that, compared with { the scalar time-fractional PDE,} the explicit prolongation formula of the infinitesimal generator for multi-dimensional fractional system (\ref{eqn-ger}) is still not quite clear and the expression of $\mu_s$ also becomes more complex. Therefore, we start with the prolongation formula { involving Riemann-Liouville fractional derivative} and then show the symmetry structure and the determining conditions of Lie symmetries for system (\ref{eqn-ger}).}
%
%$\F_t=(\partial \F_1/\partial t,\dots,\partial \F_q/\partial t)$ and $\F_{x_i}=(\partial \F_1/\partial x_i,\dots,$ $\partial \F_q/\partial x_i), i=1,2,\dots,p$, %meaning that each element in $\F$ is acted by the first order $t$-derivatives and $\x$-derivatives respectively,%are still functions of $\u$ or its $\x$-derivatives,.
%does not include the terms only containing functions of $t$ and $\x$.

The remainder of the article is outlined as follows: In Section 2, after recall the related definitions and properties of Riemann-Liouville fractional derivative, we first give an explicit prolongation formula in multi-dimensional {fractional} case and the general form of the infinitesimal generators of Lie symmetries for system (\ref{eqn-ger}), and then present two simple conditions { to determine} the infinitesimal generators. In Section 3, we use the results to study three types of time-fractional PDEs. The last section concludes the results.
\section{Main results}
\subsection{Preliminaries}
We first review the definition and some related properties of Riemann-Liouville fractional derivative, {for details please refer to} \cite{pod-1999,ki-1993}.
\begin{define}\label{def-1}
The Riemann-Liouville fractional derivative for a continuous function $u=u(t,\text{\x})$ in $[0,b]\times\mathbb{R}$ is defined by
\begin{eqnarray}
&&\no \partial_t^\alpha u=\frac{\partial^\alpha u}{\partial t^\alpha}=
\begin{cases}
\displaystyle{\frac{1}{\Gamma(m-\alpha)}\frac{\partial^m}{\partial t^m}\int_{0}^{t} (t-\nu)^{n-\alpha-1}u(\nu,\x)d\nu},\ &0\leq m-1<\alpha<m,\vspace{0.2cm}\\
\displaystyle{\frac{\partial^m u}{\partial t^m}},\ &\alpha=m\in \mathbb{N},
\end{cases}%\\
%&& \hspace{2.1cm}=\begin{cases}
%\displaystyle{\sum_{n=0}^{\infty}\binom{\alpha}{n} \frac{t^{n-\alpha}}{\Gamma(n+1-\alpha)}\frac{\partial^n}{\partial t^n}u(t,x)},\ &n-1<\alpha<n,\\
%\displaystyle{\frac{\partial^n u}{\partial t^n}},\ &\alpha=n.
%\end{cases}
\end{eqnarray}
where the gamma function is $\Gamma(z)=\int_0^{\infty}e^{-z}t^{z-1}dt$.
\end{define}

By Definition 2.1, for a power function $t^\gamma$, we have
%\begin{eqnarray}\label{proper-1}
%%&& \partial_t^\alpha(u\,v)=\sum_{k=0}^{\infty}\binom{\alpha}{k}\partial_t^ k u\,\partial^{\alpha-k}_t v,\\
%&&\no\partial_t^\alpha (t^\gamma)=\frac{\Gamma(\gamma+1)}{\Gamma(\gamma+1-\alpha)}t^{\gamma-\alpha},~~~\gamma>\alpha-1.
%%&& D_t^\alpha [f(u)]=\frac{df}{du}\,D_t^\alpha u=D_u^\alpha f(u')^\alpha.
%\end{eqnarray}
\begin{eqnarray}\label{proper-1}
&& \partial_t^\alpha \,t^\gamma=\begin{cases}
\displaystyle{\frac{\Gamma(\gamma+1)}{\Gamma(\gamma+1-\alpha)}t^{\gamma-\alpha}},~~~&\gamma>\alpha-1,\\
0,\ &\gamma=\alpha-1.
\end{cases}
\end{eqnarray}
Then for a function $f=f(\x)$ independent of $t$, $\partial_t^\alpha (f)=f\,t^{-\alpha}/\Gamma(1-\alpha)$, thus $\partial_t^\alpha f=0$ if and only if $f=0$. Let $u=u(t,\x)$ and $v=v(t,\x)$ be two continuous functions in $[0,b]\times\mathbb{R}$ along with all its $t$-derivatives of $u$, Then the Riemann-Liouville fractional derivatives of their sum and product are listed as follows \cite{ki-1993,pod-1999}
\begin{eqnarray}\label{proper}
%&&\no \partial_t^n(\partial_t^\alpha)u=\partial_t^{n+\alpha}u,\\
&& \no\partial_t^\alpha(a\,u+b\,v)=a\ \partial_t^\alpha(u)+b\,\partial_t^\alpha(v),\\
&& \partial_t^\alpha(u\,v)=\sum_{k=0}^{\infty}\binom{\alpha}{k}\partial_t^ k u\,\partial^{\alpha-k}_t v,
%&&\partial_t^\alpha t^\gamma=\frac{\Gamma(\gamma+1)}{\Gamma(\gamma+1-\alpha)}t^{\gamma-\alpha},~~~\gamma>\alpha-1.
%&& D_t^\alpha [f(u)]=\frac{df}{du}\,D_t^\alpha u=D_u^\alpha f(u')^\alpha.
\end{eqnarray}
where $a,b$ are two constants, and the second equality is called the generalized Leibniz rule, hereinafter,
\begin{eqnarray}
&&\no\binom{\alpha}{k}=\frac{(-1)^{k-1}\alpha\Gamma(k-\alpha)}{\Gamma(1-\alpha)\,\Gamma(k+1)}.
\end{eqnarray}

%In the equality $\partial_t^m(\partial_t^\alpha\,u)=\partial_t^{m+\alpha}u=\partial_t^\alpha(\partial_t^m\,u)$, the former holds identically while the latter works if and only if $u^{(k)}=0$ with $k=0,1,\dots,m-1$.
In particular for $v=1$ and the infinite differentiable function $u=u(t,\x)$, by the generalized Leibniz rule in (\ref{proper}) the Riemann-Liouville fractional derivative can be expressed as
\begin{eqnarray}\label{for-1}
&& \partial_t^\alpha u=\begin{cases}
\displaystyle{\sum_{k=0}^{\infty}\binom{\alpha}{k} \frac{t^{k-\alpha}}{\Gamma(k+1-\alpha)}\frac{\partial^k}{\partial t^k}u},\ &0\leq m-1<\alpha<m,\vspace{0.15cm}\\
\displaystyle{\frac{\partial^m u}{\partial t^m}},\ &\alpha=m.
\end{cases}
\end{eqnarray}

{ Meanwhile, we use $D_t^\alpha$ to denote the fractional total derivative with respect to $t$ { and define it as }\cite{rosa-2017}
\begin{eqnarray} \label{tofrac}
D_t^\alpha=\sum_{k=0}^{\infty}\binom{\alpha}{k} \frac{t^{k-\alpha}}{\Gamma(k+1-\alpha)}D_t^k,
\end{eqnarray}
where $D_t$ is the total derivative in $t$ and satisfies $D_t^0(u)=u$, $D_t^{n+1}=D_t(D_t^n)$. In the case of { $u=u(t,x)$ with} two independent variables $t$ and $x$,  $D_t$ is defined as
 \begin{eqnarray}
&&\no D_t=\partial_t+u_t\partial_u+u_{xt}\partial_{u_{x}}+u_{tt}\partial_{u_t}+\dots.
\end{eqnarray}
\begin{lemma}\label{tro}
Let $u=u(t,x)$ {and} $v=v(t,x)$ be two infinite differentiable functions.  Then $D_t^\alpha$ satisfies the generalized Leibniz rule
\begin{eqnarray}\label{tofrac-lei}
&& D_t^\alpha(u\,v)=\sum_{i=0}^{\infty}\binom{\alpha}{i}D_t^i u\,D^{\alpha-i}_t v,
\end{eqnarray}
where the operator $D^{\alpha-i}_t$ is defined by replacing $\alpha$ in (\ref{tofrac}) { with} $\alpha-i$.
\end{lemma}

\emph{Proof.} Acting the fractional total derivative $D_t^\alpha$ on the product $(u\,v)$ yields
\begin{eqnarray} \label{tofrac1}
&&\no D_t^\alpha(u\,v)=\sum_{k=0}^{\infty}\binom{\alpha}{k} \frac{t^{k-\alpha}}{\Gamma(k+1-\alpha)}D_t^k(u\,v)\\
&&\no \hspace{1.45cm}=\sum_{k=0}^{\infty}\binom{\alpha}{k} \frac{t^{k-\alpha}}{\Gamma(k+1-\alpha)}\sum_{i=0}^{k}\binom{k}{i} D_t^iu\,D_t^{k-i}v\\
&&\no \hspace{1.45cm}=\sum_{i=0}^{\infty} D_t^iu\sum_{k=i}^{\infty}\binom{\alpha}{k}\binom{k}{i} \frac{t^{k-\alpha}}{\Gamma(k+1-\alpha)}\,D_t^{k-i}v\\
&&\no \hspace{1.45cm}\xlongequal{k=j+i}\sum_{i=0}^{\infty}\binom{\alpha}{i} D_t^iu\sum_{j=0}^{\infty}\binom{\alpha-i}{j}\frac{t^{j+i-\alpha}}{\Gamma(j+i+1-\alpha)}\,D_t^jv\\
&&\no \hspace{1.45cm}=\sum_{i=0}^{\infty}\binom{\alpha}{i}D_t^i u\,D^{\alpha-i}_t v,
\end{eqnarray}
where the Leibniz rule of integer-order derivative is used in the second step. The proof ends. $\hfill{} \Box$}

%The generalized chain rule for fractional derivative of a composite function $f(g(t))$ is given by \cite{pod-1999}
%\begin{eqnarray}\label{formu-1}
%&& \frac{d^\alpha f(g(t))}{d t^\alpha}=\sum_{n=0}^\infty\sum_{k=0}^n\frac{1}{n!} \binom{n}{k}[-g(t)]^k\,\partial_t^\alpha[g^{n-k}(t)]\frac{d^n f(g)}{d g^n}.
%\end{eqnarray}
%where
%\begin{eqnarray}
%&&\no U_n=\sum_{k=0}^n(-1)^k \binom{n}{k}g^k(t)\partial_t^\alpha(g^{n-k}(t)).
%\end{eqnarray}
\subsection{Prolongation formula}
 Give a one-parameter local Lie symmetry group of infinitesimal transformation
\begin{eqnarray}\label{group}
&&\no t^*=t+\epsilon \,\tau(t,\textbf{x},\textbf{u})+O(\epsilon^2),\\
&&\no x_i^*=x_i+\epsilon \,\xi_i(t,\textbf{x},\textbf{u})+O(\epsilon^2),\\
&&u_{s}^*=u_{s}+\epsilon
\,\eta_{s}(t,\textbf{x},\textbf{u})+O(\epsilon^2),
\end{eqnarray}
with the group parameter $\epsilon$, which is completely characterized by the infinitesimal generator \cite{olv,blu,lv-1982}
\begin{align}\label{eq:rel6}
\mathcal{X}=\tau(t,\x,\u)\partial_t+\xi_i(t,\x,\u)\partial_{x_i}+\eta_{s}(t,\x,\u)\partial_{u_{s}},
\end{align}
where the infinitesimals $\tau=\tau(t,\x,\u),\xi_i=\xi_i(t,\x,\u)$ and $\eta_s=\eta_s(t,\x,\u)$ are determined by
\begin{eqnarray}
&&\no \tau=\frac{dt^*}{d \epsilon}_{|\epsilon=0},~~~\xi_i=\frac{dx_i^*}{d \epsilon}_{|\epsilon=0},~~~\eta_s=\frac{du_s^*}{d \epsilon}_{|\epsilon=0},
\end{eqnarray}
$i=1,\dots,p,\, s=1,\dots,q$. Thus finding Lie symmetry group (\ref{group}) is equivalent to determine the infinitesimal generator (\ref{eq:rel6}), i.e. $\tau,\,\xi_i$ and $\eta_s$. Note that here and in the rest of the paper we assume that $0<\alpha<1$ and the summation convention for repeated indices is used unless otherwise noted.
%\begin{theorem}
%(Lie theorem for FPDE) Let $\mathscr{E}=0$ be a $(1+1)$-dimensional FPDE, if $$\mbox{Pr}^{(\alpha,n)}X(\mathscr{E})|_{\mathscr{E}=0},$$ then
%\end{theorem}

Following the Lie invariance criterion for time-fractional PDEs \cite{rosa-2017,by-1998,gaz-2007}, we {find that if} system (\ref{eqn-ger}) is admitted by the Lie symmetry group (\ref{group}), then the corresponding infinitesimal generator $\mathcal{X}$ in (\ref{eq:rel6}) satisfies two conditions
\begin{align}\label{eq:rel7}
\mbox{Pr}^{(\alpha,\,k)}\mathcal{X}\left(\partial_t^\alpha \u-\F-\H\right)|_{\{\partial_t^\alpha \u-\F-\H=0\}}=0
\end{align}
and $\tau(t,\x,\u)|_{t=0}=0$, where $|_{\Delta}$ means { that the evaluations work under the condition} $\Delta$, $\mbox{Pr}^{(\alpha,\,l)}\mathcal{X}$ denotes the prolongation of the infinitesimal generator $\mathcal{X}$ in (\ref{eq:rel6}) and {is defined by}
\begin{eqnarray}\label{pro-1}
\text{Pr}^{(\alpha,\,k)}\mathcal{X}=\mathcal{X}+\eta_s^\alpha\frac{\partial}{\partial(\partial_t^\alpha u_s)}+\sum_\theta\eta_s^{\theta}(t,\textbf{x},\textbf{u}^{(k)}) \frac{\partial}{\partial u_s^{\theta}},
\end{eqnarray}
with the second summation being over all $\theta=\{\theta_1,\dots,\theta_p\}\in Z_{+}^p$ and $1\leq|\theta|\leq k$. The coefficient functions $\eta_s^{\theta}(t,\textbf{x},\textbf{u}^{(k)})$ are given by the formula \cite{olv,blu}
\begin{eqnarray}\label{pro-2}
\eta_s^{\theta}(t,\textbf{x},\textbf{u}^{(k)})=D_\theta\left(\eta_s-\tau \partial_tu_s-\sum_{i=0}^p\xi_i\,u_s^i\right)+\tau \partial_tu_s^{\theta}+\sum_{i=0}^p\xi_i\,u_s^{\theta,\,i},
\end{eqnarray}
where $u_s^i=\partial u_s/\partial x_i$ and $u_s^{\theta,\,i}=\partial u_s^{\theta}/\partial x_i$, $D_\theta=D_{\theta_1}\dots D_{\theta_p}$ is the $|\theta|$-order total derivative operator with respect to \textbf{x} and $D_{\theta_i}=D^{\theta_i}_{x_i}$.
 The symbols $D_{x_i}$ denote the total derivatives with respect to $x_i$,
 \begin{eqnarray}
&&\no D_{x_i}=\partial_{x_i}+u_{x_i}\partial_u+u_{x_ix_j}\partial_{u_{x_j}}+u_{x_it}\partial_{u_t}+\dots.
\end{eqnarray}
Now we give an explicit expression of $\eta_s^\alpha$ in the case of $p$ dependent variables and $q$ independent variables.
\begin{lemma}\label{lemma-1}
The {coefficient function $\eta_s^\alpha$ in (\ref{pro-1})} related to Riemann-Liouville time-fractional derivative is expressed by
\begin{eqnarray}\label{for-eta-1}
&& \eta_s^\alpha=D_t^\alpha (\eta_s-\tau \partial_{t}u_s-\xi_i \partial_{x_i}u_s)+\tau \partial_t^{\alpha+1} (u_s)+\xi_i\,\partial_t^\alpha (\partial_{x_i}u_s),
%&&\hspace{0.5cm}=D_t^\alpha (\eta_s)-D_t^\alpha(\tau \partial_{t}u_s)-D_t^\alpha(\xi_i \partial_{x_i}u_s)+\tau\partial_t(\partial_t^{\alpha} u_s)+\xi_i\partial_{x_i}(\partial_t^\alpha u_s),
%D_t^\alpha (\eta_s-\tau \partial_{t}u_s-\xi_i \partial_{x_i}u_s)+\tau\partial_t(\partial_t^{\alpha} u_s)+\xi_i\,\partial_{x_i}(\partial_t^\alpha u_s)\\
%&&\hspace{0.5cm},
\end{eqnarray}
where  $D_t^\alpha$ is the total fractional derivative with respect to $t$ and defined by (\ref{tofrac}).
 %In order to give an explicit expression of $\eta_s^\alpha$, we focus on the
%fractional total derivative $D_t^\alpha$.
\end{lemma}

\emph{Proof.} Extending the transformation group (\ref{group}) to fractional derivative $\partial_t^{\alpha}u_s$  as well as using the series expression (\ref{for-1}) of fractional derivative, we find
\begin{eqnarray}\label{lem1-eq1}
&& \no\hspace{-0.5cm}\eta_s^\alpha=\frac{d}{d \epsilon}\left[\frac{\partial^\alpha}{\partial (t^*)^\alpha}u_s^*(t^*,\x^*)\right]_{|\epsilon=0}\\
&&\no =\frac{d}{d \epsilon}\left[\sum_{n=0}^{\infty}\binom{\alpha}{n} \frac{(t^*)^{n-\alpha}}{\Gamma(n+1-\alpha)}\,\frac{\partial^n}{\partial (t^*)^n}u_s^*(t^*,\x^*)\right]_{|\epsilon=0}\\
%&&\no =\sum_{n=0}^{\infty}\binom{\alpha}{n} \frac{d}{d \epsilon}\left[\frac{(t+\epsilon \,\tau)^{n-\alpha}}{\Gamma(n+1-\alpha)}\,\frac{\partial^n}{\partial (t^*)^n}u_s^*(t^*,\x^*,\u^*)\right]_{|\epsilon=0}\\
%&&\no =\sum_{n=0}^{\infty}\binom{\alpha}{n} \frac{1}{\Gamma(n+1-\alpha)}\left[(n-\alpha)t^{n-\alpha-1}\tau\frac{\partial^nu_s}{\partial t^n}+t^{n-\alpha}\frac{\partial^n \eta}{\partial t^n}\right]\\
&& =\sum_{n=0}^{\infty}\binom{\alpha}{n} \frac{t^{n-\alpha-1}\tau }{\Gamma(n-\alpha)}\frac{\partial^nu_s}{\partial t^n}+\sum_{n=0}^{\infty}\binom{\alpha}{n} \frac{t^{n-\alpha}}{\Gamma(n+1-\alpha)}\eta_s^{(t,\,n)},
\end{eqnarray}
where $\eta_s^{(t,\,n)}$ {is determined by }
\begin{eqnarray}
&&\no \eta_s^{(t,\,n)}=D_t^n\left(\eta_s-\tau \frac{\partial u_s}{\partial t}-\xi_i\,\frac{\partial u_s}{\partial x_i}\right)+\tau \frac{\partial^{n+1} u_s}{\partial t^{n+1}}+\xi_i\frac{\partial^{n+1} u_s}{\partial t^n\partial x_i}.
\end{eqnarray}
Then substituting it into (\ref{lem1-eq1}), we have
\begin{eqnarray}\label{lem1-eq2}
&& \no\hspace{-1cm}\eta_s^\alpha=\tau \sum_{n=0}^{\infty}\binom{\alpha}{n} \frac{t^{n-\alpha-1}}{\Gamma(n-\alpha)}\frac{\partial^nu_s}{\partial t^n}+\sum_{n=0}^{\infty}\binom{\alpha}{n} \frac{t^{n-\alpha}}{\Gamma(n+1-\alpha)}D_t^n\left(\eta_s-\tau \frac{\partial u_s}{\partial t}-\xi_i\,\frac{\partial u_s}{\partial x_i}\right)\\
&&\no +\tau \sum_{n=0}^{\infty}\binom{\alpha}{n} \frac{t^{n-\alpha}}{\Gamma(n+1-\alpha)}\frac{\partial^{n+1} u_s}{\partial t^{n+1}}+\xi_i\sum_{n=0}^{\infty}\binom{\alpha}{n} \frac{t^{n-\alpha}}{\Gamma(n+1-\alpha)}\frac{\partial^{n}}{\partial t^n}\left(\frac{\partial u_s}{\partial x_i}\right)\\
&&\no\hspace{-0.5cm} =D_t^\alpha\left(\eta_s-\tau \frac{\partial u_s}{\partial t}-\xi_i\,\frac{\partial u_s}{\partial x_i}\right)+\xi_i\frac{\partial^\alpha }{\partial t^\alpha}\left(\frac{\partial u_s}{\partial x_i}\right)\\
&&\no+\frac{\tau \,t^{-1-\alpha}u_s}{\Gamma(-\alpha)}+\tau\sum_{n=1}^{\infty}\binom{\alpha}{n} \frac{t^{n-\alpha-1} }{\Gamma(n-\alpha)}\frac{\partial^nu_s}{\partial t^n}+\tau \sum_{n=1}^{\infty}\binom{\alpha}{n-1} \frac{t^{n-1-\alpha}}{\Gamma(n-\alpha)}\frac{\partial^{n} u_s}{\partial t^{n}}\\
&&\no\hspace{-0.5cm} =D_t^\alpha\left(\eta_s-\tau \frac{\partial u_s}{\partial t}-\xi_i\,\frac{\partial u_s}{\partial x_i}\right)+\xi_i\frac{\partial^\alpha }{\partial t^\alpha}\left(\frac{\partial u_s}{\partial x_i}\right)+\tau\sum_{n=0}^{\infty}\binom{\alpha+1}{n} \frac{t^{n-\alpha-1} }{\Gamma(n-\alpha)}\frac{\partial^nu_s}{\partial t^n}\\
&&\no \hspace{-0.5cm}=D_t^\alpha\left(\eta_s-\tau \frac{\partial u_s}{\partial t}-\xi_i\,\frac{\partial u_s}{\partial x_i}\right)+\xi_i\frac{\partial^\alpha }{\partial t^\alpha}\left(\frac{\partial u_s}{\partial x_i}\right)+\tau \frac{\partial^{\alpha+1}u_s}{\partial t^{\alpha+1}},
\end{eqnarray}
where property (\ref{for-1}) is used repeatedly.
The proof ends. $\hfill{} \Box$

Note that the integer-order prolongation formula in \cite{olv} is immediately recovered in the limit case $\alpha\rightarrow1$ while the $\alpha$th-order one in $(1+N)$-dimensional case in \cite{rosa-2017} is obtained by choosing $s=1$.
\begin{lemma}\label{lemma-2}
An explicit expression of time-fractional total derivative $D_t^\alpha$ in terms of { Riemann-Liouville} fractional derivative is given by
\begin{eqnarray}\label{fra-1}
&&\no D_t^\alpha (\eta_s)=\frac{\partial^\alpha \eta_s}{\partial t^\alpha} +\sum_{i=1}^q\left[\frac{\partial\eta_s}{\partial u_i }\frac{\partial^\alpha u_i}{\partial t^\alpha}-u_i\frac{\partial^\alpha }{\partial t^\alpha} \left(\frac{ \partial\eta_s}{\partial u_i}\right)\right]\\
&&\hspace{2.65cm}+\sum_{i=1}^q\sum_{n=1}^{\infty}\binom{\alpha}{n}\frac{\partial^n}{\partial t^n}\left(\frac{\partial \eta_s}{\partial u_i}\right)\partial_t^{\alpha-n}(u_i)+\mu_s,
\end{eqnarray}
where
\begin{eqnarray}\label{mu}
&&\mu_s=\sum_{n=2}^{\infty}\binom{\alpha}{n} \frac{t^{n-\alpha}}{\Gamma(n+1-\alpha)} \sum_{m_1+\dots+m_q=2}^n\binom{n}{m_1}\prod_{j=1}^{q-1}\binom{n-\sum_{b=1}^jm_b}{m_{j+1}} \\
&&\no\hspace{1cm}\times\sum_{k_1=0}^{m_1}\dots\sum_{k_q=0}^{m_q}\prod_{i=1}^q\left[\,\sum_{r_i=0}^{k_i}\frac{1}{k_i!}\binom{k_i}{r_i}(-u_i)^{r_i}\frac{\partial ^{m_i}}{\partial t^{m_i}}\left(u_i^{k_i-r_i}\right)\right]\frac{\partial^{m_0} }{\partial t^{m_0}}\left(\frac{\partial^k \eta_s}{\partial
u_1^{k_1}\partial u_2^{k_2}\dots
\partial u_q^{k_{q}}}\right)
\end{eqnarray}
with the indexes satisfying $k_1+\dots+k_q=k\geq2$ and $m_0+m_1+m_2+\dots+m_q=n$.
\end{lemma}
%\begin{align}
%%&\nonumber\binom{\alpha}{n}=\frac{(-1)^{n-1}\alpha\Gamma(n-\alpha)}{\Gamma(1-\alpha)\,\Gamma(n+1)},
%\no\mu=\sum_{n=2}^{\infty}\sum_{m=2}^{n}\sum_{k=2}^{m}\sum_{r=0}^{k-1}\binom{\alpha}{n}&\binom{n}{m}\binom{k}{r}
%\frac{1}{k!}\frac{t^{n-\alpha}(-u)^r}{\Gamma(n+1-\alpha)} \frac{\partial^m}{\partial t^m}(u^{k-r})\frac{\partial^{n-m}}{\partial t^{n-m}}\left(\frac{\partial^{k}\eta}{\partial u^k}\right).
%\end{align}
\emph{Proof.} By means of the  generalized Leibniz rule of fractional derivative in (\ref{proper}) and the generalized chain rule of integer-order derivative for a composite function, we obtain
\begin{eqnarray}\label{total1}
&&\no \hspace{-0.6cm}D_t^\alpha(\eta_s)=\sum_{n=0}^{\infty}\binom{\alpha}{n} \frac{t^{n-\alpha}}{\Gamma(n+1-\alpha)} ~D_t^n(\eta_s)\\
&&\no\hspace{0.7cm}=\sum_{n=0}^{\infty}\binom{\alpha}{n} \frac{t^{n-\alpha}}{\Gamma(n+1-\alpha)} \sum_{m_1+\dots+m_q=0}^n\binom{n}{m_1}\\
&&\no \hspace{1cm}\times\prod_{a=1}^{q-1}\binom{n-\sum_{b=1}^am_b}{m_{a+1}}\frac{\partial^n \eta_s\left(t,\x,u_1(t_1,\x),\dots,u_q(t_q,\x)\right)}{\partial t^{m_0}\partial
t_1^{m_1}\partial t_2^{m_2}\dots
\partial t_q^{m_{q}}}|_{\{t_1=t,\,\dots ,\,t_q=t\}}\\
&&\no\hspace{0.7cm}=\sum_{n=0}^{\infty}\binom{\alpha}{n} \frac{t^{n-\alpha}}{\Gamma(n+1-\alpha)} \sum_{m_1+\dots+m_q=0}^n\binom{n}{m_1}\prod_{j=1}^{q-1}\binom{n-\sum_{b=1}^j m_b}{m_{j+1}}\sum_{k_1=0}^{m_1}\\
&&\hspace{1.2cm}\dots\sum_{k_q=0}^{m_q}\prod_{i=1}^q\left[\sum_{r_i=0}^{k_i}\frac{1}{k_i!}\binom{k_i}{r_i}(-u_i)^{r_i}\frac{\partial ^{m_i}}{\partial t^{m_i}}\left(u_i^{k_i-r_i}\right)\right]\frac{\partial^{m_0} }{\partial t^{m_0}}\left(\frac{\partial^k \eta_s}{\partial
u_1^{k_1}\partial u_2^{k_2}\dots
\partial u_q^{k_{q}}}\right),
\end{eqnarray}
where $k_1+\dots+k_q=k$ and $m_0+m_1+m_2+\dots+m_q=n$.

By means of the Fa$\grave{\mbox{a}}$ di Bruno formula for the $m_i$th-order $t$-derivative of $u_i^{k_i-r_i}$ and the direct computations, we find
\begin{eqnarray}\label{exp-mus}
\sum_{r_i=0}^{k_i}\binom{k_i}{r_i}(-u_i)^{r_i}\frac{\partial ^{m_i}}{\partial t^{m_i}}\left(u_i^{k_i-r_i}\right)=\sum_{r_i=0}^{k_i} \binom{k_i}{r_i}(-1)^{r_i}\sum \binom{k_i-r_i}{a}\frac{ (m_i)!\,a!\,u_i^{k_i-a}}{a_1!\dots a_{m_i}!}\prod_{j=1}^{m_i}\left(\frac{\partial_t^j u_i}{j!}\right)^{a_j},
\end{eqnarray}
where $a_j$ are nonnegative integers, the second sum works on $a_1+\dots+a_{m_i}=a\leq (k_i-r_i)$ and $a_1+2a_2+\dots+m_ia_{m_i}=m_i$. Moreover, each term in (\ref{exp-mus}) is homogeneous in $u_i$ and its $t$-derivatives and the total degree is $k_i$.
Thus we isolate the linear terms in $\u$ and its derivatives in (\ref{total1}), which is equivalent to search for the terms with the indexes satisfying $k=k_1+\dots+k_q\leq1$. The solutions of this inequality are divided into $q+1$ cases where the first one is $k_i=0$ for each $i\in Z_+^{(q)}$, the other $q$ cases are $k_i=1$ for some $i\in Z_+^{(q)}$ and $k_j=0$ for { all other} $j(\neq i)\in  Z_+^{(q)}$, where $ Z_+^{(q)}$ denotes the set of positive integer not greater than $q$.

Case I. For each $i\in  Z_+^{(q)}$, $k_i=0$ means $k=0$ and requires $m_i=0,m_0=n$ for all $i\in  Z_+^{(q)}$, otherwise $\partial^{m_i}\big(u_i^{k_i-r_i}\big)/\partial t^{m_i}=0$. Thus this case corresponds to the term $\mbox{P}_{I}=\partial^\alpha \eta_s/\partial t^\alpha$, { where property (\ref{for-1}) is used.}

Case II. Consider the last $q$ cases. We fix $i$ and analyze $k_i=1$ and all $k_j=0$ with $j(\neq i)\in  Z_+^{(q)}$. Then similar as Case I, for each $j\neq i$, $k_j=0$ requires $m_j=0$ and then $k_i=k=1, m_0+m_i=n$ where $m_i\in [1,n]$ is an arbitrary positive integer since $m_i\geq k_i=1$. Meanwhile, the value of $r_i$ associated with $k_i$ is divided into two cases $k_i=r_i=1$ and $k_i=1,r_i=0$. In the former case, since $u_i^{k_i-r_i}=u_i^0=1$ and $m_i\geq 1$, then $\partial^{m_i}\big(u_i^{k_i-r_i}\big)/\partial t^{m_i}=0$ and thus { all terms vanish}. For the latter case, the terms with the given $i$ are collected as
\begin{eqnarray}
&&\no \mbox{P}_{II}^i=\sum_{n=0}^{\infty}\binom{\alpha}{n} \frac{t^{n-\alpha}}{\Gamma(n+1-\alpha)} \sum_{m_i=1}^n\binom{n}{m_i}\frac{\partial ^{m_i}}{\partial t^{m_i}}\left(u_i\right)\frac{\partial^{n-m_i}}{\partial t^{n-m_i}}\left(\frac{\partial\eta_s}{\partial
u_i}\right)\\
&&\no \hspace{0.7cm}=\sum_{n=0}^{\infty}\binom{\alpha}{n} \frac{t^{n-\alpha}}{\Gamma(n+1-\alpha)}\left[ \frac{\partial^n}{\partial t^n}\left(u_i\frac{\partial\eta_s}{\partial u_i }\right)-u_i\frac{\partial^n}{\partial t^n} \left(\frac{ \partial\eta_s}{\partial u_i}\right)\right]\\
&&\no \hspace{0.7cm}=\frac{\partial^\alpha}{\partial t^\alpha}\left(u_i\frac{\partial\eta_s}{\partial u_i }\right)-u_i\frac{\partial^\alpha }{\partial t^\alpha} \left(\frac{ \partial\eta_s}{\partial u_i}\right)\\
&&\no\hspace{0.7cm}=\frac{\partial\eta_s}{\partial u_i }\frac{\partial^\alpha u_i}{\partial t^\alpha}+\sum_{n=1}^{\infty}\binom{\alpha}{n}\frac{\partial^n}{\partial t^n}\left(\frac{\partial \eta_s}{\partial u_i}\right)\partial_t^{\alpha-n}(u_i)-u_i\frac{\partial^\alpha }{\partial t^\alpha} \left(\frac{ \partial\eta_s}{\partial u_i}\right),
\end{eqnarray}
where property (\ref{for-1}) and the generalized Leibniz rule in (\ref{proper}) are used in the third and last steps respectively.

Therefore, separating the term $\mbox{P}_{I}$ in Case I and the terms $\sum_{i=1}^q\mbox{P}_{II}^i$ in Case II from $D_t^\alpha(\eta_s)$, we obtain the explicit expression of $D_t^\alpha(\eta_s)$ in (\ref{fra-1}).
It completes the proof. $\hfill{} \Box$

%Observe that from the proof of Lemma \ref{lemma-2}, $\mu_s$ is nonlinear in $\u$ and its $\x$-derivatives.
By Lemmas \ref{lemma-1} and \ref{lemma-2}, we give an explicit expression of $\eta_s^\alpha$.
\begin{theorem} \label{th-eta}
An explicit expression of { the coefficient function} $\eta_s^\alpha$ is given by
\begin{eqnarray}\label{for-eta-22}
&&\no\eta_s^\alpha= \frac{\partial^\alpha \eta_s}{\partial t^\alpha} +\sum_{i=1}^q\left[\frac{\partial\eta_s}{\partial u_i }\frac{\partial^\alpha u_i}{\partial t^\alpha}-u_i\frac{\partial^\alpha }{\partial t^\alpha} \left(\frac{ \partial\eta_s}{\partial u_i}\right)\right]+\sum_{i=1}^q\sum_{k=1}^{\infty}\binom{\alpha}{k}\frac{\partial^k}{\partial t^k}\left(\frac{\partial \eta_s}{\partial u_i}\right)\partial_t^{\alpha-k}(u_i)\\
&&\hspace{1cm}-\sum_{k=0}^{\infty}\binom{\alpha}{k+1}D_t^{k+1} (\tau)\,\partial^{\alpha-k}_t (u_s)-\sum_{k=1}^{\infty}\binom{\alpha}{k}D_t^ k (\xi_i)\,\partial^{\alpha-k}_t \left(\frac{\partial u_s}{\partial x_i }\right)+\mu_s,
\end{eqnarray}
where $\mu_s$ is given by (\ref{mu}).
\end{theorem}

\emph{Proof.} Using the generalized Leibniz rule  { in Lemma \ref{tro} for $D_t^\alpha$}, we have
\begin{eqnarray}\label{fra-2}
&&\no  D_t^\alpha(\xi_i \partial_{x_i}u_s)=\xi_i\,D^{\alpha}_t (\partial_{x_i}u_s)+\sum_{k=1}^{\infty}\binom{\alpha}{k}D_t^ k (\xi_i)\,D^{\alpha-k}_t (\partial_{x_i}u_s),\\
&&\no D_t^\alpha(\tau \partial_{t}u_s)=\sum_{k=0}^{\infty}\binom{\alpha}{k}D_t^{k} (\tau)\,D^{\alpha-k+1}_t (u_s)\\
&&\no\hspace{1.85cm}=\tau D_t^{\alpha+1} (u_s)+\sum_{k=1}^{\infty}\binom{\alpha}{k}D_t^{k} (\tau)\,D^{\alpha-k+1}_t (u_s)\\
&& \hspace{1.85cm}=\tau D_t^{\alpha+1} (u_s)+\sum_{k=0}^{\infty}\binom{\alpha}{k+1}D_t^{k+1} (\tau)\,D^{\alpha-k}_t (u_s).
\end{eqnarray}

Then inserting (\ref{fra-1}) and  (\ref{fra-2}) into (\ref{for-eta-1}) yields (\ref{for-eta-22})
since $D_t^{\alpha-k}(u_i)=\partial_t^{\alpha-k}(u_i),D^{\alpha-k}_t (\partial_{x_i}u_s)=\partial^{\alpha-k}_t (\partial_{x_i}u_s)$.
It completes the proof. $\hfill{} \Box$

{ The expression of $\eta_s^\alpha$ in (\ref{for-eta-22}) includes the previous prolongation} formulas in \cite{ss-2018,pra-2017,kd-2018,sg-2017} as special cases and also revises some inaccurate expressions of $\mu_s$.
In particular, for the $(1+1)$-dimensional case, i.e. $x_1=x,u_1=u,\eta_s=\eta$, then $\alpha$th-order prolongation (\ref{for-eta-22}) becomes \cite{zhang-2019}
\begin{align}
\no \eta^\alpha=&\nonumber\partial_t^\alpha\eta+\big[\eta_u-\alpha D_t(\tau)\big]\partial_t^\alpha u-u\,\partial_t^\alpha (\eta_u)+\mu- \sum_{k=1}^{\infty}\binom{\alpha}{k}D_t^k(\xi)\partial_t^{\alpha-k}(u_x)
 \\\no &+\sum_{k=1}^{\infty}
 \left[\binom{\alpha}{k}\partial_t^k(\eta_u)-\binom{\alpha}{k+1}D_t^{k+1}(\tau)\right]
 \partial_t^{\alpha-k}(u),
\end{align}
where %$D^\alpha_t$ is the total fractional derivative with respect to $t$, and
\begin{align}
%&\nonumber\binom{\alpha}{n}=\frac{(-1)^{n-1}\alpha\Gamma(n-\alpha)}{\Gamma(1-\alpha)\,\Gamma(n+1)},
\no\mu=\sum_{n=2}^{\infty}\sum_{m=2}^{n}\sum_{k=2}^{m}\sum_{r=0}^{k-1}\binom{\alpha}{n}&\binom{n}{m}\binom{k}{r}
\frac{1}{k!}\frac{t^{n-\alpha}(-u)^r}{\Gamma(n+1-\alpha)} \frac{\partial^m}{\partial t^m}(u^{k-r})\frac{\partial^{n-m}}{\partial t^{n-m}}\left(\frac{\partial^{k}\eta}{\partial u^k}\right).
\end{align}

{Therefore, by means of coefficients formulas (\ref{pro-2}) and (\ref{for-eta-22}), we get an explicit expression of the prolongation $\mbox{Pr}^{(\alpha,\,k)}\mathcal{X}$ in (\ref{pro-1}) for multi-dimensional fractional case.} Then the procedure for searching the infinitesimals $\tau,\xi_i$ and $\eta_s$ in (\ref{eq:rel6}) is similar as the one of integer-order PDEs. Thus substituting the formula (\ref{pro-1}) into condition (\ref{eq:rel7}) and then annihilating to zero first the coefficients of time-fractional integrals and derivatives of $\u$ and followed by the coefficients of integer-order $\x$-derivatives of $\u$, we obtain an over-determined system for $\tau,\,\xi_i$ and $\eta_s$ which includes integer-order derivative and fractional { integral} and derivative operations. Then solving the {determining} system together with the condition $\tau(t,\x,\u)|_{t=0}=0$ gives the infinitesimal generator (\ref{eq:rel6}).
\subsection{Symmetry structure}
It is well-known that finding solutions of the determining system of Lie symmetry is a rather labor-consuming task for integer-order PDEs \cite{blu,olv,lv-1982}, not to mention herein the system involving fractional { integral and} derivative operations. Thus in this section, we analyze the symmetry structure of system (\ref{eqn-ger}) and  { show that the infinitesimal generator $\mathcal{X}$ in (\ref{eq:rel6}) has a simple and unified expression}. Such a  scenario of knowing the structure of infinitesimal generator in advance will simplify the symmetry determining system largely. Note that in the procedure of searching for Lie symmetries, all time-fractional integrals and derivatives arising in the prolongation formulas are considered as independent variables.

%We show that the infinitesimal generator (\ref{eq:rel6}) admitted by Eq.(\ref{eqn-ger}) takes a simplified form, which definitely drops off the whole computational cost.
\begin{lemma}\label{lemm-a}
Let $\mu_s$ be given in (\ref{mu}). A necessary and sufficient condition for $\mu_s=0$ is that $\eta_s=\eta_s(t,\x,\u)$ is linear in $u_i$.
\end{lemma}

\emph{Proof}. By the expression of $\mu_s$ in (\ref{mu}), if $\eta_s$ is linear in $u_i$ then $\partial^k \eta_s/\partial u_1^{k_1}\dots\partial u_q^{k_q}=0$ with $k=k_1+\dots+k_q\geq2$ and thus $\mu_s=0$. The sufficiency holds. Next we prove necessity, i.e. prove $\partial^2 \eta_s/\partial u_i^2=\partial^2 \eta_s/\partial u_i \partial u_l=0\,(i\neq l)$. By the proof of Lemma \ref{lemma-2}, we find that all terms in $\mu_s$ are nonlinear in $u_i$ and their $t$-derivatives.

Consider the term $(\partial_tu_i)^2$ which occurs uniquely for $m_i=2,\,m_j=0\,(j\neq i)$, $i,j=1,\dots,q$.
Thus we separate the case $m_i=2$ from $\mu_s$ and rewrite $\mu_s$ as the following form
\begin{eqnarray}\label{mu-sepa}
&&\no\hspace{-0.5cm}\mu_s=\sum_{n=2}^{\infty}\binom{\alpha}{n}\binom{n}{2} \frac{t^{n-\alpha}}{\Gamma(n+1-\alpha)} \sum_{r_i=0}^{2}\left[\frac{1}{2!}\binom{2}{r_i}(-u_i)^{r_i}\frac{\partial ^{2}}{\partial t^{2}}\left(u_i^{2-r_i}\right)\right]\frac{\partial^{n-2} }{\partial t^{n-2}}\left(\frac{\partial^2 \eta_s}{\partial
u_i^2}\right)+\mbox{remainder}\\%
%&&\no\hspace{0.6cm}+\sum_{n=2}^{\infty}\binom{\alpha}{n} \frac{t^{n-\alpha}}{\Gamma(n+1-\alpha)} \sum_{\substack{m_1+\dots+m_q=2\\ m_i\neq 2}}^n\binom{n}{m_1}\prod_{j=1}^{q-1}\binom{n-\sum_{b=1}^j m_b}{m_{j+1}} \\
%&&\no\hspace{0.6cm}\times\sum_{k_1=0}^{m_1}\dots\sum_{k_q=0}^{m_q}\prod_{i=1}^q\left[\sum_{r_i=0}^{k_i}\frac{1}{k_i!}\binom{k_i}{r_i}(-u_i)^{r_i}\frac{\partial ^{m_i}}{\partial t^{m_i}}\left(u_i^{k_i-r_i}\right)\right]\frac{\partial^{m_0} }{\partial t^{m_0}}\left(\frac{\partial^k \eta_s}{\partial
%u_1^{k_1}\partial u_2^{k_2}\dots
%\partial u_q^{k_{q}}}\right)\\
&&=\left[\,\sum_{n=2}^{\infty}\binom{\alpha}{n} \binom{n}{2} \frac{t^{n-\alpha}}{\Gamma(n+1-\alpha)}\frac{\partial^{n-2} }{\partial t^{n-2}}\left(\frac{\partial^2 \eta_s}{\partial
u_i^2}\right)\right](\partial_tu_i)^2+\mbox{remainder}.%\\
%&&\hspace{0.6cm}+\sum_{n=2}^{\infty}\binom{\alpha}{n} \frac{t^{n-\alpha}}{\Gamma(n+1-\alpha)} \sum_{\substack{m_1+\dots+m_q=2\\ m_i\neq 2}}^n\binom{n}{m_1}\prod_{j=1}^{q-1}\binom{n-\sum_{b=1}^jm_b}{m_{j+1}} \\
%&&\no\hspace{0.6cm}\times\sum_{k_1=0}^{m_1}\dots\sum_{k_q=0}^{m_q}\prod_{i=1}^q\left[\sum_{r_i=0}^{k_i}\frac{1}{k_i!}\binom{k_i}{r_i}(-u_i)^{r_i}\frac{\partial ^{m_i}}{\partial t^{m_i}}\left(u_i^{k_i-r_i}\right)\right]\frac{\partial^{m_0} }{\partial t^{m_0}}\left(\frac{\partial^k \eta_s}{\partial
%u_1^{k_1}\partial u_2^{k_2}\dots
%\partial u_q^{k_{q}}}\right),
\end{eqnarray}
%where the indexes $k_1+\dots+k_q=k\geq2$ and $m_0+m_1+m_2+\dots+m_q=n$.
For the coefficient of $(\partial_tu_i)^2$ in (\ref{mu-sepa}), let $\lambda=n-2$, then it becomes
\begin{eqnarray}\label{sum-1-coe}
&&\no\sum_{\lambda=0}^{\infty}\binom{\alpha}{\lambda+2}\binom{\lambda+2}{2}
\frac{t^{\lambda+2-\alpha}}{\Gamma(\lambda+3-\alpha)}\frac{\partial ^{\lambda}}{\partial t^{\lambda}}\left(\frac{\partial^2 \eta_s}{\partial
u_i^2}\right)\\
&&\no=\frac{1}{2}\alpha(\alpha-1)\sum_{\lambda=0}^{\infty}\binom{\alpha-2}{\lambda}
\frac{t^{\lambda+2-\alpha}}{\Gamma(\lambda+3-\alpha)}\frac{\partial ^{\lambda}}{\partial t^{\lambda}}\left(\frac{\partial^2 \eta_s}{\partial
u_i^2}\right)\\
&&=\frac{1}{2}\alpha(\alpha-1) \partial_t^{\alpha-2}\left(\frac{\partial^2 \eta_s}{\partial
u_i^2}\right),
\end{eqnarray}
where property (\ref{for-1}) is used. By the uniqueness of $(\partial_tu_i)^2$, annihilating its coefficient to zero yields $\partial_t^{\alpha-2}(\partial^2 \eta_s/\partial u_i^2)=0$ which {means} $\partial^2 \eta_s/\partial u_i^2=C(\x,\u)\,t^{\alpha-3}$ with an undetermined function $C(\x,\u)$. Then further {splitting} the case $k=k_i=2$ from $\mu_s$, which implies for $j\neq i$, $k_j=m_j=0$, we rewrite  $\mu_s$ as the form %, if some $k_i\geq 3$ the term $\partial^{k_i} \eta_s/\partial u_i^{k_i}=0$. Thus $\mu_s$ is further simplified to
\begin{eqnarray}\label{mu-sepa-new1}
&&\no\hspace{-0.5cm}\mu_s=C(\x,\u)\Bigg\{\sum_{n=3}^{\infty}\binom{\alpha}{n} \frac{t^{n-\alpha}}{\Gamma(n+1-\alpha)} \sum_{m_i=3}^n\binom{n}{m_i} \\
&&\no\hspace{1.8cm}\times\sum_{r_i=0}^{2}\left[\frac{1}{2!}\binom{2}{r_i}(-u_i)^{r_i}\frac{\partial ^{m_i}}{\partial t^{m_i}}\left(u_i^{2-r_i}\right)\right]\frac{\partial^{m_0} }{\partial t^{m_0}}\left(\, t^{\alpha-3}\right)\Bigg\}+\mbox{remainder}\\
&&\no=C(\x,\u)\Bigg[\sum_{n=3}^{\infty}\binom{\alpha}{n} \frac{t^{n-\alpha}}{\Gamma(n+1-\alpha)} \sum_{m_i=3}^n\binom{n}{m_i}\frac{\partial^{n-m_i} }{\partial t^{n-m_i}}\left(\, t^{\alpha-3}\right)\Bigg] \\
&&\hspace{1.8cm}\times\left[\sum_{\kappa=1}^{m_i-1}\binom{m_i-1}{\kappa}\partial_t^\kappa u_i\,\partial_t^{m_i-\kappa} u_i\right]+\mbox{remainder}.
\end{eqnarray}

Then by interchanging the order of summations and adopting the technique used in (\ref{sum-1-coe}), we further arrange (\ref{mu-sepa-new1}) as the form
\begin{eqnarray}\label{mu-sepa-new}
&&\no\hspace{-0.5cm}\mu_s=C(\x,\u)\Bigg[\sum_{j=3}^{\infty}\sum_{n=j}^{\infty}\binom{\alpha}{n} \frac{t^{n-\alpha}}{\Gamma(n+1-\alpha)}\binom{n}{j}\frac{\partial^{n-j} }{\partial t^{n-j}}\left(\, t^{\alpha-3}\right)\Bigg] \\
&&\no\hspace{1.8cm}\times\left[\sum_{\kappa=1}^{m_i-1}\binom{m_i-1}{\kappa}\partial_t^\kappa u_i\,\partial_t^{m_i-\kappa} u_i\right]+\mbox{remainder}\\
&&\no =C(\x,\u)\Bigg[\sum_{j=3}^{\infty}\binom{\alpha}{j} \sum_{k=0}^{\infty}\binom{\alpha-j}{k}\frac{t^{k+j-\alpha}}{\Gamma(k+j+1-\alpha)}\frac{\partial^k }{\partial t^k}\left(\, t^{\alpha-3}\right)\Bigg] \\
&&\no\hspace{1.8cm}\times\left[\sum_{\kappa=1}^{m_i-1}\binom{m_i-1}{\kappa}\partial_t^\kappa u_i\,\partial_t^{m_i-\kappa} u_i\right]+\mbox{remainder}\\
&&=C(\x,\u)\sum_{j=3}^{\infty}\binom{\alpha}{j} \partial_t^{\alpha-j}\left( t^{\alpha-3}\right)\left[\sum_{\kappa=1}^{m_i-1}\binom{m_i-1}{\kappa}\partial_t^\kappa u_i\,\partial_t^{m_i-\kappa} u_i\right]+\mbox{remainder}.
\end{eqnarray}

By the uniqueness of nonlinear terms $\partial_t^\kappa u_i\,\partial_t^{m_i-\kappa} u_i$ with $m_i\neq 2$ in $\mu_s$, we obtain the coefficients $C(\x,\u)=0$ and then $\partial^2 \eta_s/\partial u_i^2=0$. Thus in $\mu_s$, the terms $\partial^k \eta_s/\partial u_i^k\,(k\geq2)$ { and} their derivatives vanish while the remaining terms take the form $\partial^k \eta_s/\partial u_1^{k_1}\partial u_2^{k_2}\dots\partial u_q^{k_{q}}$ with nonnegative integer $k_i\leq 1 \,,i=1,\dots,q$.

Next consider the cross derivative terms $\partial^2 \eta_s/\partial u_i \partial u_l\, (i\neq l)$. It corresponds to $k=2, k_i=k_l=1$ and $k_d=0$  with $d(\neq i,l)=1,\dots ,q$, thus $m_i\geq1, m_l\geq1$ and $m_d=0$, otherwise $\mu_s$ vanishes identically. Then one has
\begin{eqnarray}
\no\left[\,\sum_{r_i=0}^{1}\binom{1}{r_i}(-u_i)^{r_i}\frac{\partial ^{m_i}}{\partial t^{m_i}}\left(u_i^{1-r_i}\right)\right]\left[\,\sum_{r_l=0}^{1}\binom{1}{r_l}(-u_l)^{r_l}\frac{\partial ^{m_l}}{\partial t^{m_l}}\left(u_l^{1-r_l}\right)\right]=\partial_t^{m_i}u_i\,\partial_t^{m_l}u_l,
\end{eqnarray}
which are unique in $\mu_s$.

We first consider $m_i=m_l=1$, i.e. the term $\partial_t u_i\partial_t u_l$. Then $m_0=n-2$ and the coefficient of $\partial_t u_i\partial_t u_l$  is
\begin{eqnarray}\label{sum-4}
\sum_{n=2}^{\infty}\binom{\alpha}{n} \frac{n(n-1) t^{n-\alpha}}{\Gamma(n+1-\alpha)}\,\partial_t^{n-2}\left(\frac{\partial^{2}\eta_s}{\partial u_i\partial u_l}\right)=\alpha(\alpha-1) \partial_t^{\alpha-2}\left(\frac{\partial^2 \eta_s}{\partial u_i \partial u_l}\right),
\end{eqnarray}
where the technique adopted in (\ref{sum-1-coe}) is used again. {Then by solving $\partial^{\alpha-2}_t\left(\partial^2 \eta_s/\partial u_i \partial u_l\right)=0$, we get} $\partial^2 \eta_s/\partial u_i \partial u_l=B(\x,\u)\,t^{\alpha-3}$ with an undetermined function $B(\x,\u)$.

Secondly, consider the case $m_i=m_l=2$ which corresponds to the term $\partial_t^2 u_i\partial_t^2 u_l$. Then $m_0=n-4$, and the coefficient of $\partial_t^2 u_i\partial_t^2 u_l$ is
\begin{eqnarray}
&&\no \frac{1}{4}\alpha(\alpha-1)(\alpha-2)(\alpha-3)\,\partial_t^{\alpha-4}\left(\frac{\partial^2 \eta_s}{\partial u_i \partial u_l}\right)=\frac{1}{4}\alpha(\alpha-1)(\alpha-2)(\alpha-3)B(\x,\u)\,\partial_t^{\alpha-4}\left(t^{\alpha-3}\right)\\
&&\hspace{7cm} {=\no \frac{1}{4}(\alpha-3)\Gamma(\alpha+1)tB(\x,\u)}.
 \end{eqnarray}
Thus by the uniqueness of $\partial_t^2 u_i\partial_t^2 u_l$ we obtain $B(\x,\u)=0$ and then  { for all $i\neq j$,} $\partial^2 \eta_s/\partial u_i \partial u_l=0$.
It completes the proof. $\hfill{} \Box$

\begin{theorem}\label{th-1}
If the infinitesimal generator $\mathcal{X}$ given in (\ref{eq:rel6}) leaves system (\ref{eqn-ger}) invariant, then $\mathcal{X}$ must take the form
\begin{equation}\label{ger-form}
 \mathcal{X}=(\chi_2t^2+\chi_1 t)\partial_t+\xi_i(\text{\x})\partial_{x_i}+\eta_s(t,\x,u)\partial_{u_s},
\end{equation}
where $\chi_1$ and $\chi_2$ are arbitrary constants, $\eta_s=\eta_s(t,\x,\u)$  is given by
\begin{equation}\label{gamma}
\eta_s=\big[\,g_s(\x)+\gamma(2\chi_2t+\chi_1)\big]u_s+\sum_{i\neq s}f_i(\x)u_i+h_s(t,\x),
\end{equation}
where $g_s(\x),f_i(\x)$ and $h_s(t,\x)$ are undetermined smooth functions of their arguments respectively, and the constant $\gamma$ satisfies
\begin{equation}\label{gamma-1}
\gamma=\begin{cases}
0,\ &\chi_2= 0,\\
\displaystyle{\frac{1}{2}(\alpha-1)},\ &\chi_2\neq 0.
\end{cases}
\end{equation}
\end{theorem}

\emph{Proof}. We show the theorem by analyzing the structure of Eq.(\ref{eq:rel7}) on the space $(\partial_t^{\alpha}\u^{(k)},$ $\partial_t^{\alpha-1}\u^{(k)}, \dots)$. Thus expanding Eq.(\ref{eq:rel7}) on the solution space of system (\ref{eqn-ger}) yields
\begin{equation}\label{eq:rel12}
\eta_s^\alpha-\tau(t,\x,\u) (\F_t+\H_t)-\xi_i(t,\x,\u) (\F_{x_i}+\H_{x_i})-\sum_{j=1}^q\sum_\theta\eta_j^{\theta}(t,\textbf{x},\textbf{u}^{(k)}) \partial_{u_j^{\theta}}(\F)=0,% ~~~\mbox{when}~ \partial_t^\alpha u-F=0.
\end{equation}
where $\eta_s^\alpha$ take the form %In Eq.(\ref{eq:rel12}) and next equations, we use the notations $u_0=u,\eta^{(i)}=\eta$.
\begin{eqnarray}\label{for-eta-2}
&&\no\eta_s^\alpha= \frac{\partial^\alpha \eta_s}{\partial t^\alpha} +\sum_{i=1}^q\left[\frac{\partial\eta_s}{\partial u_i }(\F_i+\H_i)-u_i\,\partial_t^\alpha \left(\frac{ \partial\eta_s}{\partial u_i}\right)\right]\\
&&\no \hspace{1cm}+\sum_{i=1}^q\sum_{k=1}^{\infty}\binom{\alpha}{k}\partial_t^k\left(\frac{\partial \eta_s}{\partial u_i}\right)\partial_t^{\alpha-k}(u_i)-\alpha D_t(\tau)(\F_s+\H_s)\\
&&\hspace{1cm}-\sum_{k=1}^{\infty}\binom{\alpha}{k+1}D_t^{k+1} (\tau)\,\partial^{\alpha-k}_t (u_s)-\sum_{k=1}^{\infty}\binom{\alpha}{k}D_t^k (\xi_i)\,\partial^{\alpha-k}_t \left(\frac{\partial u_s}{\partial x_i }\right)+\mu_s,
\end{eqnarray}
and $\eta_j^{\theta}(t,\textbf{x},\textbf{u}^{(k)}) $ is given by formula (\ref{pro-2}).

Then substituting (\ref{pro-2}) and (\ref{for-eta-2}) into Eq.(\ref{eq:rel12}) and vanishing the coefficients of $\partial^{\alpha-k}_t(\partial u_s/\partial x_i)$, one obtains
\begin{equation}\label{eqn-sy-1}
\binom{\alpha}{k}D_t^k (\xi_i)=0,~~k=1,2,\dots.
\end{equation}
Since {equations} (\ref{eqn-sy-1}) work for $k=1,2,\dots$, then for $k=1$ we have
\begin{equation}
\no D_t(\xi_i)=\frac{\partial\xi_i}{\partial t}+\sum_{j=1}^q\frac{\partial\xi_i}{\partial u_j}\frac{\partial u_j}{\partial t}=0,
\end{equation}
which implies $\partial\xi_i/\partial t=\partial\xi_i/\partial u_j=0$, i.e. $\xi_i=\xi_i(\x)$. Then equations (\ref{eqn-sy-1}) with $k\geq 2$ hold identically.

Next consider $\tau=\tau(t,\x,\u)$. We claim that in system (\ref{eqn-ger}), for each $x_i$, there exists at least one $u_j$ such that $\partial u_j/\partial x_i$ or its higher order $\x$-derivatives {occur}. If not, $x_i$ can be regarded as a parameter variable and system (\ref{eqn-ger}) involves $p-1$ independent variables, which is contradictory. Thus we assume { that such type of derivative in $\F$ has the maximal order $|\vartheta|$ and takes} the form $u_j^\vartheta$ with $\vartheta=(\vartheta_1,\dots,\vartheta_p)\in Z_{+}^p$ satisfying $1\leq|\vartheta|\leq k$. Observe that in Eq.(\ref{eq:rel12}) the derivative $\partial_t u_j^{\widehat{\vartheta}}$ uniquely exists in $\eta_j^{\vartheta}(t,\textbf{x},\textbf{u}^{(k)})$ given by (\ref{pro-2}), where $\widehat{\vartheta}=(\vartheta_1,\dots,\vartheta_{i-1},\vartheta_i-1,\vartheta_{i+1},\dots,\vartheta_p)$. More precisely, the term $\partial_t u_j^{\widehat{\vartheta}}$  appears uniquely in $ D_\vartheta (\tau \partial_t u_j)$ and its coefficient is $\partial_{u_j^{\vartheta}}(\F)D_{x_i}\tau$ which should be vanished, i.e. \begin{eqnarray}
&&\no D_{x_i}\tau=\frac{\partial\tau}{\partial x_i}+\sum_{j=1}^q\frac{\partial u_j}{\partial x_i}\frac{\partial \tau}{\partial u_j}=0,
\end{eqnarray}
 since $\partial_{u_j^{\vartheta}}(\F)\neq0$ by the claim. Note that {by direct computations} the term $\tau \partial_tu_j^{\vartheta}$ in (\ref{pro-2}) vanishes identically. Thus we obtain $\partial_{x_i}\tau=\partial_{u_j}\tau=0$, i.e. $\tau=\tau(t)$. Then the prolonged formula (\ref{pro-2}) is simplified to
\begin{eqnarray}\label{pro-2-si}
&& \eta_j^{\theta}(t,\textbf{x},\textbf{u}^{(k)})=D_\theta\left(\eta_j-\sum_{i=0}^p\xi_i\,u_j^i\right)+\sum_{i=0}^p\xi_i\,u_j^{\theta,i}.
\end{eqnarray}

With the above simplifications, condition (\ref{eq:rel12}) becomes
\begin{eqnarray}\label{eq:rel12-sim}
&&\no\partial_t^\alpha \eta_s +\sum_{i=1}^q\left[\frac{\partial\eta_s}{\partial u_i }(\F_i+\H_i)-u_i\partial_t^\alpha\left(\frac{ \partial\eta_s}{\partial u_i}\right)\right]+\sum_{i=1}^q\sum_{k=1}^{\infty}\binom{\alpha}{k}\partial_t^k\left(\frac{\partial \eta_s}{\partial u_i}\right)\partial_t^{\alpha-k}(u_i)\\
&&\no \hspace{1cm}-\sum_{k=1}^{\infty}\binom{\alpha}{k+1}D_t^{k+1} (\tau)\,\partial^{\alpha-k}_t (u_s)+\mu_s-\alpha D_t(\tau)(\F_s+\H_s)\\
&&\hspace{1cm}-\tau(t) (\F_t+\H_t)-\xi_i(\x) (\F_{x_i}+\H_{x_i})-\sum_{j=1}^q\sum_\theta\eta_j^{\theta}(t,\textbf{x},\textbf{u}^{(k)}) \partial_{u_j^{\theta}}(\F)=0,
\end{eqnarray}
where $\eta_j^{\theta}(t,\textbf{x},\textbf{u}^{(k)}) $ is given by (\ref{pro-2-si}).

We now turn to consider $\eta_s=\eta_s(t,\x,\u)$. Observe that integer-order $t$-derivatives of $u_i$ in Eq.(\ref{eq:rel12-sim}) uniquely occur in $\mu_s$ while $\partial_t^{\alpha-k}(u_i)$ with integers $k>1$ are fractional { integrals of} $u_i$. Thus by considering the structure of $\mu_s$ and Lemma \ref{lemm-a}, we obtain
%
%
%Then together with $\partial^2 \eta_s/\partial
%u_i^2=0$, we obtain $\mu_s=0$ and assume
\begin{eqnarray}\label{eta-1}
\eta_s=\sum_{i=1}^qr_i(t,\x)u_i+h_s(t,\x).
 \end{eqnarray}
where $r_i(t,\x)$ and $h_s(t,\x)$ are undetermined functions.%the symbol $\u\setminus u_s=\{u_1,\dots,u_{s-1},u_{s+1},\dots,u_q\}$ denotes the set by deleting $u_s$ from $\u$.

Next { we further separate} Eq.(\ref{eq:rel12-sim}) with respect to time-fractional integrals and derivatives of $u_i$  and get
\begin{subequations}\label{eq20}
\begin{equation}\label{eqn-sy-2}
\partial_t^{\alpha-k}u_s:~~\binom{\alpha}{k}\partial_t^k\left(\frac{ \partial\eta_s}{\partial u_s}\right)-\binom{\alpha}{k+1}D_t^{k+1}(\tau)=0,
\end{equation}
\begin{equation}\label{eqn-sy-3}
\partial_t^{\alpha-k}u_i:~~ \binom{\alpha}{k}\partial_t^k\left(\frac{ \partial\eta_s}{\partial u_i}\right)=0, ~~i(\neq s)=1,2,\dots,q,
\end{equation}
\begin{equation}\label{eqn-sy-4}
\begin{aligned}
 \mbox{Else}:  \frac{\partial^\alpha \eta_s}{\partial t^\alpha} +\sum_{i=1}^q\left[\frac{\partial\eta_s}{\partial u_i }(\F_i+\H_i)-u_i\frac{\partial^\alpha }{\partial t^\alpha} \left(\frac{ \partial\eta_s}{\partial u_i}\right)\right]-\alpha D_t(\tau)(\F_s+\H_s)\\
-\tau\left(\frac{\partial\F}{\partial t}+\frac{\partial\H}{\partial t}\right)-\xi_i \left(\frac{\partial\F}{\partial x_i}+\frac{\partial\H}{\partial x_i}\right)-\sum_{j=1}^q\sum_\theta\eta_j^{\theta}\partial_{u_j^{\theta}}(\F)=0,
\end{aligned}
\end{equation}
\end{subequations}
which hold for $k=1,2,\cdots$. %Solving system (\ref{eq20}) will get the required expression of $X$ in (\ref{eq:rel6}).

Since equations (\ref{eqn-sy-3}) holds for $k=1,2,\dots$, then for $k=1$, { together with} (\ref{eta-1}) one has $\partial^2\eta_s/\partial t\partial u_i=\partial r_i/\partial t=0$ with $i\neq s$, thus
\begin{eqnarray}\label{eta-fin}
\eta_s=r_s(t,\x)u_s+\sum_{i\neq s}f_i(\x)u_i+h_s(t,\x),
 \end{eqnarray}
 where $f_i(\x)$ are undetermined functions.

With such results equations (\ref{eqn-sy-3}) with integers $k\geq 2$ hold identically. Next we consider equations (\ref{eqn-sy-2}) to further {find} explicit expressions of $\tau$ and $\eta_s$. Since $\tau=\tau(t)$, we isolate the case $k=1$ {from equations (\ref{eqn-sy-2}) and divide them} into two parts
\begin{eqnarray}\label{eqn-sy-2-1}
&&\no \frac{\partial^2\eta_s}{\partial t\partial u_s}-\frac{1}{2}(\alpha-1)\tau''=0,~~~k=1,\\
&& \frac{\partial^{k}}{\partial t^{k}}\left(\frac{\partial\eta_s}{\partial u_s}\right)-\frac{\alpha-k}{k+1}\frac{\partial^{k+1}}{\partial t^{k+1}}(\tau)=\frac{\partial^{k-1}}{\partial t^{k-1}}\left(\frac{\partial^2\eta_s}{\partial t\partial u_s}-\frac{\alpha-k}{k+1}\tau''\right)=0,~~~k\geq2.
\end{eqnarray}
We use the partial derivative on $\tau$ in order to write the second part uniformly. Then inserting the first equation into the second ones gives
\begin{eqnarray}\label{eq-au}
&& \left(\frac{\alpha-1}{2}-\frac{\alpha-k}{k+1}\right)\frac{d^{k-1}}{d t^{k-1}}(\tau'')=0,~~~k=2,3,\cdots.
\end{eqnarray}

In particular, for $k=2$, equations (\ref{eq-au}) give $d^{3}\tau(t)/dt^{3}=0$. Together with the condition $\tau(t)|_{t=0}=0$ we obtain $\tau=\chi_1t+\chi_2t^2$, where $\chi_1$ and $\chi_2$ are two arbitrary constants.
By means of the above results, solving the first equation in system (\ref{eqn-sy-2-1}) yields two different cases.

(I). One is $\tau''=0$, i.e. $\chi_2=0$, then $\partial^2\eta_s/\partial t\partial u_s=0$. By considering (\ref{eta-fin}), we obtain
$\tau=\chi_1t,~\partial r_s(t,\x)/\partial t=0$. Thus
\begin{equation}
\no \eta_s=g_s(\x) u_s+\sum_{i\neq s}f_i(\x)u_i+h_s(t,\x).
\end{equation}

(II). The other is $\tau''\neq0$, i.e. $\chi_2\neq0$. From (\ref{eta-fin}) and the first equation in system (\ref{eqn-sy-2-1}) we obtain $r_s(t,\x)=(\alpha-1)\tau'/2+g_s(\x)$, then
\begin{equation}
\no \eta_s=\left[g_s(\x)+\frac{1}{2}(\alpha-1)\tau'\right]u_s+\sum_{i\neq s}f_i(\x)u_i+h_s(t,\x).
\end{equation}

{Finally,} we collect the above two expressions of $\eta_s$ as a unified form $\eta_s=\left[\,g_s(\x)+\gamma\tau'\,\right]u_s+\sum_{i\neq s}f_i(\x)u_i+h_s(t,\x)$ where $\gamma=(\alpha-1)/2$ for $\chi_2\neq0$ and $\gamma=0$ for $\chi_2=0$, the functions $g_s(\x),f_i(\x)$ and $h_s(t,\x)$  are determined by Eq.(\ref{eqn-sy-4}). It completes the proof. $\hfill{} \Box$
\subsection{Determining conditions}
By means of the symmetry structure of system (\ref{eqn-ger}), we show that Lie symmetries of system (\ref{eqn-ger}) are determined by two {elegant} conditions which provide a possibility to use the known computer programs of integer-order PDEs to solve the symmetry determining equations of multi-dimensional time-fractional PDEs.

\begin{theorem}\label{cor-1}
Following the above notations, Lie symmetries of system (\ref{eqn-ger}) are completely determined by %a linear fractional PDE and an integer-order PDE
\begin{eqnarray}\label{eqn-sys-det-cor}
&&\no \frac{\partial^\alpha h_s(t,\x)}{\partial t^\alpha}+\sum_{i=1}^q \frac{\partial\eta_s}{\partial u_i }\H_i-\alpha D_t(\tau)\H_s-\tau\frac{\partial\H_s}{\partial t}-\xi_i \frac{\partial\H_s}{\partial x_i}-\sum_{j=1}^q\sum_{u_j^{\theta}\in J_s}D_{\theta}(h_j(t,\x))\partial_{u_j^{\theta}}(\F_s)=0,\\
&&\no \sum_{i=1}^q\frac{\partial\eta_s}{\partial u_i }\F_i-\alpha D_t(\tau)\F_s-\tau\frac{\partial\F_s}{\partial t}-\xi_i \frac{\partial\F_s}{\partial x_i}\\
&& \hspace{1.75cm}-\sum_{j=1}^q\sum_{u_j^{\theta}\in J_s  }\left[\eta_j^{\theta}-D_{\theta}(h_j(t,\x))\right]\partial_{u_j^{\theta}}(\F_s)-\sum_{j=1}^q\sum_{u_j^{\theta}\in I_s\setminus J_s }\eta_j^{\theta}\partial_{u_j^{\theta}}(\F_s)=0,
\end{eqnarray}
where $I_s=\{\mbox{all terms in} \,\,\F_s\}$, $J_s=\{\mbox{all terms in}\, \F_s \mbox{which are linear in}\, u_j^\theta\}$ and $I_s\setminus J_s=\{$the terms contained in $\,I_s$ not in \,$J_s\}$.
\end{theorem}

\emph{Proof}. By the proof of Theorem \ref{th-1}, Eq.(\ref{eqn-sy-4}) becomes
\begin{eqnarray}\label{eqn-sys-det}
&&\no \frac{\partial^\alpha \eta_s}{\partial t^\alpha} +\sum_{i=1}^q\left[\frac{\partial\eta_s}{\partial u_i }(\F_i+\H_i)-u_i\frac{\partial^\alpha }{\partial t^\alpha} \left(\frac{ \partial\eta_s}{\partial u_i}\right)\right]-\alpha D_t(\tau)(\F_s+\H_s)\\
&& \hspace{2.2cm}-\tau\left(\frac{\partial\F_s}{\partial t}+\frac{\partial\H_s}{\partial t}\right)-\xi_i \left(\frac{\partial\F_s}{\partial x_i}+\frac{\partial\H_s}{\partial x_i}\right)-\sum_{j=1}^q\sum_\theta\eta_j^{\theta}\partial_{u_j^{\theta}}(\F_s)=0,
\end{eqnarray}
where $\tau=\chi_2t^2+\chi_1 t$, $\xi_i=\xi_i(x)$ and $\eta_s$ is given by (\ref{gamma}). Moreover,  { equations (\ref{eqn-sy-1}), (\ref{eqn-sy-2}) and (\ref{eqn-sy-3})} hold identically with the given $\tau,\,\xi_i$ and $\eta_s$ in Theorem \ref{th-1}, thus Lie symmetries of system (\ref{eqn-ger}) are uniquely determined by Eq.(\ref{eqn-sys-det}).

On the space $(t,\x,\u)$, one has
\begin{eqnarray}
%&&\no\frac{ \partial\eta_s}{\partial u_i}=\\
&&\no \frac{\partial^\alpha \eta_s}{\partial t^\alpha} -\sum_{i=1}^qu_i\frac{\partial^\alpha }{\partial t^\alpha} \left(\frac{ \partial\eta_s}{\partial u_i}\right)=\frac{\partial^\alpha }{\partial t^\alpha}\left(\eta_s-\sum_{i=1}^qu_i \frac{ \partial\eta_s}{\partial u_i}\right)=\frac{\partial^\alpha}{\partial t^\alpha} h_s(t,\x).
\end{eqnarray}

Let $I_s=\{\mbox{all terms in} \,\F_s\}$ and $J_s=\{\mbox{all terms in}\, \F_s \mbox{which are linear in}\, u_j^\theta\}$, the set $I_s\setminus J_s=\{\mbox{The terms contained in} \,I_s~\mbox{not in }\,J_s\}$.
By considering whether the terms involve $\u$ and its $\x$-derivatives or not, we separate Eq.(\ref{eqn-sys-det}) into two parts given in (\ref{eqn-sys-det-cor}). The proof ends.$\hfill{} \Box$

Theorem \ref{cor-1} shows that the symmetry determining equations of system (\ref{eqn-ger}) can be divided into two parts: a system of integer-order PDEs in $\tau,\,\xi_i$ and $\eta_s$ and a system of linear time-fractional PDEs in $h_s(t,\x)$. Moreover, for the most of time-fractional PDEs, the second {condition} in (\ref{eqn-sys-det-cor}) ``almost" completely determines the admitted Lie symmetries while the first one either holds automatically or is used to check the final results.
\begin{cor}\label{cor-3}
If $\H(t,\x)=\textbf{0}\triangleq (0,0,\dots,0)_{1\times q}$ and $J_s$ is empty with $s=1,\dots,q$, then Lie symmetries admitted by system (\ref{eqn-ger}) are uniquely determined by the second equation in (\ref{eqn-sys-det-cor}).
\end{cor}

\emph{Proof}. If $\H(t,\x)=\textbf{0}$ and $J_s$ is empty, the first equation {in (\ref{eqn-sys-det-cor})} becomes $\partial_t^\alpha h_s(t,\x)=0$ which implies $h_s(t,\x)=C(\x)t^{\alpha-1}$ with an undetermined function $C(\x)$.  Then the Lie symmetries are totally determined by the second condition in (\ref{eqn-sys-det-cor}). The proof ends. $\hfill{} \Box$

%\subsection{An algorithm}
Following the above theoretical preparations, we formulate {the procedure} of finding Lie symmetries of system (\ref{eqn-ger}) as the following three steps:

\textbf{Step 1}. Assume system (\ref{eqn-ger}) is admitted by the infinitesimal generator (\ref{eq:rel6}), then by Theorem \ref{th-1}, the infinitesimals $\tau,\,\xi_i$ and $\eta_j$ are directly assumed to be the explicit {forms} (\ref{ger-form}).

\textbf{Step 2}. Finding the two conditions to determine the Lie symmetry. By Theorem \ref{cor-1}, first write down the expressions $\H_i,\F_i$ and the sets $I_i,J_i,I_i\setminus J_i$, then obtain the two determining conditions given by (\ref{eqn-sys-det-cor}).

\textbf{Step 3}. Further separation of the second condition in system (\ref{eqn-sys-det-cor}) with respect to $\u$ and its $\x$-derivatives to get the symmetry determining system about $\tau,\,\xi_i$ and $\eta_j$, then together with the first condition, solve the system to get the {infinitesimal} generator (\ref{eq:rel6}).

It should be mentioned that in \textbf{Step 3} the separation of the second condition in system (\ref{eqn-sys-det-cor}) is more deeper than the one for integer-order PDEs where the former one is divided about $\u$ and its $\x$-derivatives because the prescribed infinitesimals $\tau,\,\xi_i$ and $\eta_j$ in Theorem \ref{th-1} are independent of $\u$ and its $\x$-derivatives, while the latter one is done only with respect to $\x$-derivatives of $\u$. Thus the separation {in \textbf{Step 3} will generate a more simplified symmetry determining system.}

\section{Three examples}
We consider three examples to illustrate the efficiencies and applications of our results. { In the first subsection} we will adopt two methods to look for Lie symmetries of the time-fractional generalized Zakharov-Kuznetsov equation (\ref{zk}) in order to show the efficiencies of our results while in next two subsections we directly use our method for the other two examples.
\subsection{Time-fractional generalized Zakharov-Kuznetsov equation}
The first example is the time-fractional generalized Zakharov-Kuznetsov equation
\begin{eqnarray}\label{zk}
&& \partial_t^\alpha u+u^{\rho} u_x+u_{xxx}+u_{xyy}=0,
\end{eqnarray}
where $\rho$ is a nonzero constant and $u=u(t,x,y)$. Eq.(\ref{zk}) with $\alpha=\rho=1$ is the Zakharov-Kuznetsov equation which describes weakly nonlinear ion-acoustic wave in a strongly magnetized lossless plasma in two dimensions \cite{zk-1974}.

We assume that Eq.(\ref{zk}) is admitted by a one-parameter local Lie symmetry group with the infinitesimal generator
\begin{eqnarray}
&&\no X=\tau \partial_t+\xi \partial_x+ \psi\partial_y+\eta \partial_u,
\end{eqnarray}
where $\xi, \tau,\psi$ and $\eta$ are smooth functions of $t,x,y$ and $u$ respectively.
\subsubsection{The original method}
The Lie infinitesimal criterion for Eq.(\ref{zk}) gives
\begin{eqnarray}\label{pro-zk}
&& \text{Pr}^{(\alpha,\,3)}X\left(\partial_t^\alpha u+u^{\rho} u_x+u_{xxx}+u_{xyy}\right)|_{\{(\ref{zk})\}}=0,
\end{eqnarray}
where $\text{Pr}^{(\alpha,\,3)}$ {is given by} (\ref{pro-1}) with $k=3$. Specifically, expanding condition (\ref{pro-zk}) yields
\begin{eqnarray}\label{pro-zk1}
&& \eta^\alpha_t+n u^{n-1} u_x\eta+u^{\rho}\eta^{x}+\eta^{xxx}+\eta^{xyy}=0,
\end{eqnarray}
where $\eta^\alpha_t$ is formulated by (\ref{for-eta-22}) while $\eta^x$ and $\eta^{xxx},\eta^{xxy}$ are expressed by (\ref{pro-2}).

First assume $\mu=0$ in $\eta^\alpha_t$. Then inserting (\ref{for-eta-22}) into condition (\ref{pro-zk1}) and separating it with respect to different time-fractional { integrals and} derivatives of $u$, {one has}
\begin{eqnarray}\label{pro-zk2}
&&\no \binom{\alpha}{k}D_t^ k (\xi)=\binom{\alpha}{k}D_t^ k (\psi)=0,~~~k=1,2,\dots,\\
&&\no \binom{\alpha}{k}\frac{\partial^k}{\partial t^k}\left(\frac{\partial \eta}{\partial u}\right)-\binom{\alpha}{k+1}D_t^{k+1} (\tau)=0,~~~k=1,2,\dots,\\
&&\no \frac{\partial^\alpha \eta}{\partial t^\alpha}-\left[\frac{\partial\eta}{\partial u}-\alpha D_t (\tau)\right](u^{\rho} u_x+u_{xxx}+u_{xyy})-u\frac{\partial^\alpha }{\partial t^\alpha} \left(\frac{ \partial\eta}{\partial u}\right)\\
&&\hspace{5cm}+\rho u^{\rho-1} u_x\eta+u^{\rho}\eta^{x}+\eta^{xxx}+\eta^{xyy}=0.
\end{eqnarray}

Then separating the last equation in (\ref{pro-zk2}) with respect to $t$- and $x$-derivatives of $u$ generates the system of symmetry determining equations which {contains a huge number of equations and is very} complicated since $\tau,\xi,\psi$ and $\eta$ are functions of $t,x,y,u$. Thus it is obvious that the known symmetry structure {in advance}, such as $\tau$ independent of $x,y,u$ and $\eta$ linear in $u$, will definitely drop off the complexity of computations.

\subsubsection{Our method}
By Theorem \ref{th-1}, we directly assume
\begin{equation}
\no\tau=\chi_2t^2+\chi_1 t,\,~\xi=\xi(x,y),~\psi=\psi(x,y),~\eta=\big[g(x,y)+\gamma(2\chi_2t+\chi_1)\big]u+h(t,x,y),
\end{equation}
where $\xi,\,\psi$ and $g=g(x,y),\,h=h(t,x,y)$ are smooth undetermined functions.
Then by Theorem \ref{cor-1}, $\H_1=0,~\F_1=-u^{\rho} u_x-u_{xxx}-u_{xyy}$, $I_1=\{u^{\rho} u_x,u_{xxx},u_{xyy}\},~J_1=\{u_{xxx},u_{xyy}\}, ~I_1\setminus J_1=\{u^{\rho} u_x\}$, then the two conditions given in (\ref{eqn-sys-det-cor}) become
\begin{subequations}
\begin{equation}\label{sys-zk-1}
\frac{\partial^\alpha h}{\partial t^\alpha}+h_{xxx}+h_{xyy}=0,
\end{equation}
\begin{eqnarray}\label{sys-zk-2}
&&\no \big[g(x,y)+(\gamma-\alpha)(2\chi_2t+\chi_1)\big]\left(-u^{\rho} u_x-u_{xxx}-u_{xyy}\right)\\
&& \hspace{2cm}+(\eta^{xxx}-h_{xxx})+(\eta^{xyy}-h_{xyy})+\rho\,u^{\rho-1}u_x\eta+u^{\rho}\eta^x=0,
\end{eqnarray}
\end{subequations}
where $h=h(t,x,y)$, $\eta^x$ and $\eta^{xxx},\eta^{xxy}$ are given by
 \begin{eqnarray}%\label{sys-zk-pro}
&&\no \eta^x=D_x(\eta-\xi u_x-\psi u_y)+\xi u_{xx}+\psi u_{xy},\\
&&\no \eta^{xxx}=D_x^3(\eta-\xi u_x-\psi u_y)+\xi u_{xxxx}+\psi u_{xxxy},\\
&&\no \eta^{xyy}=D_xD_y^2(\eta-\xi u_x-\psi u_y)+\xi u_{xxyy}+\psi u_{xyyy},
\end{eqnarray}
which are more simpler than the usual ones {(\ref{pro-2})} since they do not contain function $\tau$.

Substituting them into Eq.(\ref{sys-zk-2}) and separating it with respect to different powers of $u$ and its derivatives, we obtain
\begin{eqnarray}\label{35-b}
&&\no \xi_y=\psi_x=g_x=g_y=h_x=\rho h=0,\\
&&\no \alpha  (2 t \chi_2+\chi_1)-3 \xi_x=0,\\
&&\no 2 \psi_y+\xi_x-\alpha  (2 \chi_2\, t +\chi_1)=0,\\
&& \rho g -\xi_x+(\alpha +\gamma  \rho) (2\chi_2 \,t +\chi_1)=0,
\end{eqnarray}
which are integer-order linear PDEs and very easy to be solved. Solving the system gives $\chi_2=\gamma=h=0$ and
\begin{eqnarray} \label{sol-zk}
&& \tau=\chi_1t,~~\xi=\frac{1}{3}\alpha \chi_1 x+c_1, ~~\psi=\frac{1}{3}\alpha \chi_1 y+c_2, ~~\eta=-\frac{2}{3\rho}\alpha \chi_1 u.
\end{eqnarray}

Observe that system (\ref{35-b}) obtained by separating Eq.(\ref{sys-zk-2}) with respect to $u$ together with $x$- and $y$-derivatives of $u$ completely determines the Lie symmetries of Eq.(\ref{zk}) and solutions (\ref{sol-zk}) automatically satisfy Eq.(\ref{sys-zk-1}) since $h=0$.

The direct role of the infinitesimal operator is to reduce the PDEs into lower-dimensional PDEs. Before the performance, we recall the two variables Erd\'{e}lyi-Kober fractional differential operator
\begin{align}\label{EK-op}
\left(P^{\epsilon,\ \alpha}_{\delta,\ \sigma}U\right)(\omega,\theta)=&\prod_{j=0}^{n-1}\left(\epsilon +j-\frac{1}{\delta}\omega\frac{\partial}{\partial \omega}-\frac{1}{\sigma}\theta\frac{\partial}{\partial \theta}\right)(K^{\epsilon+\alpha,
\ n-\alpha}_{\delta,\ \sigma} U)(\omega,\theta),\\
\nonumber&n=
\begin{cases}
\left[\alpha\right]+1,\ &\alpha \neq \mathbb{N},  \\
\alpha,\ &\alpha \in \mathbb{N},
\end{cases}
\end{align}
 in order to present an elegant expression for the reduced equations, where
\begin{eqnarray}
\no\left(K^{\epsilon,\ \alpha}_{\delta,\ \sigma} U\right)(\omega,\theta)=
\begin{cases}
\displaystyle{\frac{1}{\Gamma(\alpha)}}\int_{1}^{\infty} (u-1)^{\alpha-1}u^{-(\epsilon+\alpha)}U(\omega\, u^{\frac{1}{\delta}},\theta\, u^{\frac{1}{\sigma}})du,\ &\alpha>0,\\
U(\omega,\theta),\ &\alpha=0.
\end{cases}
\end{eqnarray}

%In particular, $\left(K^{\epsilon,\ \alpha}_{\delta,\ \delta} U\right)(\rho,\theta)=\left(K^{\epsilon,\ \alpha}_{\delta} U\right)(\rho,\theta)$.
Note that in the case of single independent variable operator (\ref{EK-op}) becomes the classical Erd\'{e}lyi-Kober fractional differential operator \cite{ki-1993}. { In what follows, we give an explicit procedure of constructing the reduced equation by one infinitesimal generator while for other two {ones} as well as the two examples below we direct list the reduced equations and {similarity} solutions without details.

\begin{prop}
By the {infinitesimal} generator $X=t \partial_t+\alpha \,x /3 \partial_x+  \alpha \,y/3\partial_y-2\alpha/(3\rho) u\partial_u$, we reduce Eq.(\ref{zk}) to
the form
\begin{eqnarray}\label{redu-1}
&& \left(P^{1-\frac{2\alpha}{3\rho}-\alpha,\ \alpha}_{\frac{3}{\alpha},\,\frac{3}{\alpha}}U\right)(z_1,z_2)=U^{\rho} \frac{\partial U}{\partial z_{1}} +\frac{\partial^{3} U}{\partial z_{1}^{3}} +\frac{\partial^{3} U}{\partial z_{1}\partial z_{2}^2},
\end{eqnarray}
where the similarity variables are $z_{1}=x\,t^{-\alpha/3}, z_{2}=y\,t^{-\alpha/3}, U(z_1,z_2)=u\,t^{2\alpha/3\rho}$.
\end{prop}

\emph{Proof.} The first step is to find the similarity variable $I=I(t,x,y,u)$ by solving the linear equation $X(I)=0$. The corresponding characteristic equations are
\begin{eqnarray}
\frac{dt}{t}=\frac{dx}{\frac{\alpha}{3} \,x}=\frac{dy}{\frac{\alpha}{3} \,y}=\frac{du}{-\frac{2\alpha}{3\rho} u},
\end{eqnarray}
which gives the similarity variables $z_{1}=x\,t^{-\alpha/3}, z_{2}=y\,t^{-\alpha/3}, U(z_1,z_2)=u\,t^{2\alpha/3\rho}$. Then {by the chain rule for integer-order derivative, we get}
\begin{align}\label{ux}
u_x=t^{-\frac{\alpha}{3}(\frac{2}{\rho}+1)}\frac{\partial U}{\partial z_{1}},~~u_{xxx}=t^{-\frac{\alpha}{3}(\frac{2}{\rho}+3)}\frac{\partial^{3} U}{\partial z_{1}^{3}},~~u_{xyy}=t^{-\frac{\alpha}{3}(\frac{2}{\rho}+3)}\frac{\partial^{3} U}{\partial z_{1}\partial z_{2}^2}.
\end{align}

Next we consider the time-fractional derivative $\partial_t^{\alpha}u$ with $0<\alpha<1$. Inserting the {above similarity variables} into the Riemann-Liouville  fractional derivative in Definition \ref{def-1}, $\partial_t^{\alpha}u$ can be written as
\begin{align}\label{16}
\frac{\partial^\alpha u}{\partial t^\alpha}=\frac{1}{\Gamma(1-\alpha)}\frac{\partial}{\partial t}\int_{0}^{t} (t-s)^{-\alpha}
s^{-\frac{2\alpha}{3\rho}}\ U\left(xs^{-\frac{\alpha}{3}},ys^{-\frac{\alpha}{3}}\right)ds.
\end{align}

Let $v=t/s,\ ds=-t/v^2dv$, then we transform {$\partial_t^{\alpha}u$ in} (\ref{16}) into the following form
\begin{equation}\label{17}
 \frac{\partial^\alpha u}{\partial t^\alpha}=\frac{\partial}{\partial t} \left[\frac{t^{1-\frac{2\alpha}{3\rho}-\alpha} }{\Gamma(1-\alpha)}\int_{1}^{\infty} (v-1)^{-\alpha}\,v^{-(2-\alpha-\frac{2\alpha}{3\rho})} U\left(xv^{\frac{\alpha}{3}},yv^{\frac{\alpha}{3}}\right)\ dv\right].
\end{equation}
Substituting {the Erd\'{e}lyi-Kober fractional differential operator in} (\ref{EK-op}) into Eq.(\ref{17}), we arrive at { a compact expression}
\begin{align}\label{18}
\frac{\partial^\alpha u}{\partial t^\alpha}=\frac{\partial}{\partial t}\left[t^{1-\frac{2\alpha}{3\rho}-\alpha}\left(K^{1-\frac{2\alpha}{3\rho},\ 1-\alpha}_{\frac{3}{\alpha},\,\frac{3}{\alpha}}U\right)(z_1,z_2)\right].
\end{align}

Observe that
\begin{align}
\nonumber t\frac{\partial}{\partial t}U(z_1,z_2)&=-\frac{\alpha}{3}t^{-\frac{\alpha}{3}} \left[x\frac{\partial}{\partial z_1}U(z_1,z_2)+y\frac{\partial}{\partial z_2}{U(z_1,z_2)}\right]\\
\nonumber&=-\frac{\alpha}{3}z_1\frac{\partial}{\partial z_1}U(z_1,z_2)-\frac{\alpha}{3}z_2\frac{\partial}{\partial z_2}U(z_1,z_2).
\end{align}
Then we convert Eq.(\ref{18}) into the form
\begin{align}\label{ut}
\nonumber&\frac{\partial^\alpha u}{\partial t^\alpha}=t^{1-\frac{2\alpha}{3\rho}-\alpha-1}\left[1-\frac{2\alpha}{3\rho}-\alpha-\frac{\alpha}{3}\left(z_1\frac{\partial }{\partial z_1}+z_2\frac{\partial }{\partial z_2}\right)\right]\left(K^{1-\frac{2\alpha}{3\rho},\ 1-\alpha}_{\frac{3}{\alpha},\,\frac{3}{\alpha}}U\right)(z_1,z_2)\\
&\hspace{0.75cm}=t^{-\frac{2\alpha}{3\rho}-\alpha}\left(P^{1-\frac{2\alpha}{3\rho}-\alpha,\ \alpha}_{\frac{3}{\alpha},\,\frac{3}{\alpha}}U\right)(z_1,z_2).
\end{align}

Finally, inserting { expressions }(\ref{ux}) and (\ref{ut}) into Eq.(\ref{pro-zk}), we obtain the reduced equation (\ref{redu-1}). The proof ends. $\hfill{} \Box$

{Following the above procedure, we find that} the similarity variables of the infinitesimal generator $X=\partial_x$ are $z_1=t,z_2=y, U(z_1,z_2)=u(t,x,y)$.  {Then Eq.(\ref{zk}) is converted} to the form $\partial^{\alpha} U(z_1,z_2)/\partial z_1^{\alpha} =0$. {By solving it and returning to original variables we find a solution of Eq.(\ref{zk}) in the form} $u(t,x,y)=f(y)t^{\alpha-1}$ with an arbitrary function $f(y)$.

Similarly, with the {infinitesimal} generator $X=\partial_y$, Eq.(\ref{zk}) is  reduced to
\begin{eqnarray}
&&\no\frac{\partial^{\alpha}U}{\partial z_1^{\alpha}}=U^{\rho}U_{z_2}+U_{z_2z_2z_2},
\end{eqnarray}
 where $z_1=t,z_2=x,U=u(t,x,y)$.
\subsection{Time-fractional Hirota-Satsuma coupled KdV equations}
Consider the time-fractional Hirota-Satsuma coupled KdV equations
\begin{eqnarray}\label{kdv}
&&\no \partial_t^\alpha u= u u_x+vv_x+u_{xxx},\\
&& \partial_t^\alpha v=- u v_x-2 v_{xxx},
\end{eqnarray}
whose Lie symmetry analysis has been performed in \cite{pra-2017}. Here we directly use our results to find Lie symmetries of system (\ref{kdv}).

Assume that a Lie symmetry with the infinitesimal generator
\begin{eqnarray}\label{oper-kdv}
&&\no X=\tau \partial_t+\xi \partial_x+\eta \partial_u+ \phi\partial_v,
\end{eqnarray}
where $\xi,\tau,\eta$ and $\phi$ are arbitrary smooth functions of $t,x,u$ {and} $v$ respectively, leaves system  (\ref{kdv}) invariant. { Then by Theorem \ref{th-1}, we get} $\tau=\chi_2t^2+\chi_1 t,\,\xi=\xi(x)$,
\begin{equation}\label{gamma-kdv-eta}
\no\eta=\big[g_1(x)+\gamma(2\chi_2t+\chi_1)\big]u+f_2(x)v+h_1(t,x),
\end{equation}
and
\begin{equation}\label{gamma-kdv-phi}
\no\phi=\big[g_2(x)+\gamma(2\chi_2t+\chi_1)\big]v+f_1(x)u+h_2(t,x).
\end{equation}

By Theorem \ref{cor-1}, one has
\begin{eqnarray}
&&\no \H_1(t,x)=0,\F_1=u u_x+vv_x+u_{xxx},\\
&&\no I_1=\{u u_x,vv_x,u_{xxx}\},J_1=\{u_{xxx}\}, I_1\setminus J_1=\{u u_x,vv_x\};\\
&&\no \H_2(t,x)=0,\F_2=-u v_x-2 v_{xxx},\\
&&\no I_2=\{u v_x,vu_x,v_{xxx}\},J_2=\{v_{xxx}\}, I_2\setminus J_2=\{u v_x,vu_x\}.
\end{eqnarray}
Then the two conditions in (\ref{eqn-sys-det-cor}) for the first equation in system (\ref{kdv}) become
\begin{eqnarray}\label{sys-kdv1}
&& \no\frac{\partial^\alpha h_1}{\partial t^\alpha}-\beta_3\frac{\partial^3 h_1}{\partial x^3}=0,\\
&&\no \big[g_1(x)+(\gamma-\alpha)(2\chi_2t+\chi_1)\big]\F_1+f_2(x)\F_2\\
&& \hspace{2cm}-\left(\eta^{xxx}-\frac{\partial^3 h_1}{\partial x^3}\right)-\left(u_x\eta+u\eta^x\right)-\left(v_x\phi+v\phi^x\right)=0,
\end{eqnarray}
and for the second equation become
\begin{eqnarray}\label{sys-kdv2}
&&\no \frac{\partial^\alpha h_2}{\partial t^\alpha}-\delta_3\frac{\partial^3 h_2}{\partial x^3}=0,\\
&& \big[g_2(x)+(\gamma-\alpha)(2\chi_2t+\chi_1)\big]\F_2+f_1(x)\F_1+2\left(\phi^{xxx}-\frac{\partial^3 h_2}{\partial x^3}\right)+(v_x\eta+u\phi^x)=0,
\end{eqnarray}
where $h_i=h_i(t,x)$, $\eta^x,\phi^x$ and $\eta^{xxx},\phi^{xxx}$ are given by
\begin{eqnarray}\label{sys-kdv-pro}
&&\no \eta^x=D_x(\eta-\xi u_x)+\xi u_{xx},~~~\eta^{xxx}=D_x^3(\eta-\xi u_x)+\xi u_{xxxx},\\
&&\no\phi^x=D_x(\phi-\xi v_x)+\xi v_{xx}~~~\phi^{xxx}=D_x^3(\phi-\xi v_x)+\xi v_{xxxx}.
\end{eqnarray}

Inserting the above prolongations into systems (\ref{sys-kdv1}) and (\ref{sys-kdv2}) and annihilating the coefficients of different powers of $u,v$ and their derivatives to zero, we obtain
\begin{eqnarray}\label{sym-hsk}
&&\no f_1=f_2=h_1=h_2=0,g_1'=0,g_1=g_2,\\
&&\no \alpha  (2\chi_2 t +\chi_1)-3 \xi '=0,\\
&& (\alpha+\gamma)(2\chi_2 t +\chi_1)+g_1-\xi '=0.
\end{eqnarray}

Solving the system yields
\begin{eqnarray}
&&\no \tau=\chi_1t,~~\xi=\frac{1}{3}\alpha \chi_1 x+c_1, ~~\eta=-\frac{2}{3}\alpha \chi_1 u, ~~\phi=-\frac{2}{3}\alpha \chi_1 v,
\end{eqnarray}
which is the same as the results in \cite{pra-2017}. However, the symmetry determining system (\ref{sym-hsk}) is very simple and easy to be solved.

Next we use the infinitesimal operators to construct reduced equations. The {infinitesimal} generator $X=\partial_x$ reduce system  (\ref{kdv}) to the form
\begin{eqnarray}\label{reud-kdv}
&& \frac{\partial^\alpha U(\zeta)}{\partial\zeta^\alpha}=0,~~~~\frac{\partial^\alpha V(\zeta)}{\partial\zeta^\alpha}=0,
\end{eqnarray}
where the similarity variables are $\zeta=t, \, U(\zeta)=u(t,x)$ {and} $V(\zeta)=v(t,x)$. Solving system (\ref{reud-kdv}) gives one solution of system (\ref{kdv}) {in the form} $u(t,x)=C_1t^{\alpha-1},v(t,x)=C_2t^{\alpha-1}$, where $C_1$ {and} $C_2$ are integral constants.

Induced by the operator $X=t \partial_t+\alpha/3 x \partial_x-  2\alpha /3u\partial_u-2\alpha/3 v\partial_v$ and with {the} similar procedure for Eq.(\ref{redu-1}), system (\ref{kdv}) is reduced to
\begin{eqnarray}
&&\no \left(P_{\frac{3}{\alpha}}^{1-\frac{5\alpha}{3},\,\alpha}U\right)(\zeta)=U(\zeta)U'(\zeta)+V(\zeta)h'(\zeta)+U'''(\zeta),\\
&&\no\left(P_{\frac{3}{\alpha}}^{1-\frac{5\alpha}{3},\,\alpha}V\right)(\zeta)=-U(\zeta)V'(\zeta)-2V'''(\zeta),
\end{eqnarray}
where the similarity variables are $\zeta=xt^{-\alpha/3}, \, U(\zeta)=ut^{2\alpha/3}$ {and} $V(\zeta)=vt^{2\alpha/3}$.

\subsection{Time-fractional nonlinear telegraph equations}
The third example is the time-fractional nonlinear telegraph equations with variable coefficients
\begin{eqnarray}\label{tele}
&&\no \partial_t^\alpha u=v_x,\\
&& \partial_t^\alpha v=P(u)u_x+G(u),
\end{eqnarray}
where $P(u)$ and $Q(u)$ are two smooth nonzero functions of $u$ which make system (\ref{tele}) nonlinear. In what follows, we will perform a Lie symmetry classification of system (\ref{tele}) which means {first} to classify the functions $P(u)$ and $Q(u)$ making system (\ref{tele}) admit the extended symmetries and then to determine the symmetries.

{To simplify our calculations, we use} an equivalent transformation of system (\ref{tele}) given by
\begin{eqnarray}\label{equ-tran}
&&  t^*=t,x^*=\beta_3x+\beta_4,u^*=\beta_1 u,v^*=\beta_2 v,P^*(u^*)=\frac{\beta_2\beta_3}{\beta_1} P(u), G^*(u^*)=\beta_2G(u),
\end{eqnarray}
where nonzero constants $\beta_i$ satisfy $\beta_1\beta_2=\beta_3\neq0$. { Transformation (\ref{equ-tran}) maps system (\ref{tele}) into the same form. Its main role is }that in the procedure of {Lie} symmetry classification for system (\ref{tele}), scalings of ${P(u)}$ and $G(u)$ do not affect the final classified results. For example, if ${P(u)}=\beta_1u^2+\beta_2,G(u)=\beta_3 u$, we can assume ${P(u)}=u^2+\hat{\beta}_2,G(u)=u$.

Assume that a local Lie symmetry group with the infinitesimal generator
\begin{eqnarray}\label{oper-tele}
&& X=\tau \partial_t+\xi \partial_x+\eta \partial_u+ \phi\partial_v,
\end{eqnarray}
where the infinitesimals $\xi,\tau,\eta$ and $\phi$ are undetermined smooth functions of $t,x,u$ and $v$ respectively, leaves system  (\ref{tele}) invariant.  {Then} by Theorem \ref{th-1}, $\tau=\chi_2t^2+\chi_1 t,\,\xi=\xi(x)$,
\begin{equation}\label{gamma-tele-eta}
\no \eta=\big[g_1(x)+\gamma(2\chi_2t+\chi_1)\big]u+f_2(x)v+h_1(t,x),
\end{equation}
and
\begin{equation}\label{gamma-tele-phi}
\no\phi=\big[g_2(x)+\gamma(2\chi_2t+\chi_1)\big]v+f_1(x)u+h_2(t,x),
\end{equation}
where, hereinafter, $f_i(x),g_i(x)$  {and}  $h_i(t,x)$ with $i=1,2$ are undetermined functions.
Moreover, by Theorem \ref{cor-1} for system (\ref{tele}), one has $\H_1(t,x)=\H_2(t,x)=0,\F_1=v_x,\F_2=P(u)u_x+G(u)$.

We first consider the first equation in system (\ref{tele}). By Theorem \ref{cor-1}, $I_1=\{v_x\}$ and $J_1=\{v_x\}$, then the set $I_1\setminus J_1$ is empty and we obtain
\begin{eqnarray}\label{sys-tele-1}
&&\no \frac{\partial^\alpha h_1(t,x)}{\partial t^\alpha}-\frac{\partial h_2(t,x)}{\partial x}=0,\\
&& \big[g_1(x)+(\gamma-\alpha) (2\chi_2t+\chi_1)\big]v_x -\left(\phi^x-\frac{\partial h_2(t,x)}{\partial x}\right)+f_2(x)\big[P(u)u_x+G(u)\big]=0,
\end{eqnarray}
where $\phi^x$ is expressed by
\begin{eqnarray}\label{sys-tele-pro-1}
\phi^x=\big[ g_2(x)+\gamma(2\chi_2t+\chi_1)\big]v_x+g_2'(x)v+f_1'(x)u+f_1(x)u_x+\frac{\partial h_2(t,x)}{\partial x}-\xi'(x)v_x.
\end{eqnarray}

Inserting (\ref{sys-tele-pro-1}) into the second equation in system (\ref{sys-tele-1}) and separating it with respect to $v,v_x$ and $u_x$, we obtain  $g_2(x)=c_1$, and
\begin{eqnarray}\label{sys-tele-symdet-1}
&&\no f_2(x)P(u)-f_1(x)=0,\\
&&\no f_2(x)G(u)- f_1'(x)u=0,\\
&& -\alpha (2\chi_2t+\chi_1)+g_1(x)-c_1+\xi'(x)=0,
\end{eqnarray}
 {where, here and below, $c_j$ are arbitrary constants, $j=1,2,3$.}

Observe that if $f_1(x)f_2(x)\neq0$, $P(u)$ is a constant and $G(u)$ is either linear in $u$ or a constant, which is contradict with  {nonlinear system (\ref{tele}).} Thus $f_1(x)f_2(x)=0$. Since $f_2(x)P(u)=f_1(x)$ and $P(u)\neq0$, thus $f_1(x)=0$  implies $f_2(x)=0$ and vice versa, i.e. $f_1(x)=f_2(x)=0$. Then by dividing the last equation in system (\ref{sys-tele-symdet-1}) with respect to $t$ and reconsidering (\ref{sys-tele-1}), we obtain the symmetry determining equations for the first equation in system (\ref{tele})
%\begin{eqnarray}\label{sys-tele-symdet-11aa}
%&&\no \frac{\partial^\alpha h_1(t,x)}{\partial t^\alpha}-\frac{\partial h_2(t,x)}{\partial x}=0,\\
%&&\no (\gamma-\gamma-\alpha)\chi_2=0,\\
%&& (\gamma-\gamma-\alpha)\chi_1+g_1(x)-c_1+\xi'(x)=0.
%\end{eqnarray}
%
%We claim that $\chi_2=0$. If $\chi_2\neq0$, by (\ref{gamma-1}) one has $\gamma=\gamma=(\alpha-1)/2$, then the second equation in system  (\ref{sys-tele-symdet-11aa}) becomes $-\alpha \chi_2=0$, which does not hold since $\alpha\neq0$. Thus system (\ref{sys-tele-symdet-11aa}) becomes
\begin{eqnarray}\label{sys-tele-symdet-11}
&&\no \frac{\partial^\alpha h_1(t,x)}{\partial t^\alpha}-\frac{\partial h_2(t,x)}{\partial x}=0,\\
&& \chi_2=-\alpha\chi_1+g_1(x)-c_1+\xi'(x)=0.
\end{eqnarray}
Then by the above analysis the infinitesimals in (\ref{oper-tele}) are simplified to
\begin{equation}\label{gamma-tele-eta-new}
\no\tau=\chi_1 t,\ \xi=\xi(x),\  \eta=g_1(x)u+h_1(t,x),\ \phi=c_1v+h_2(t,x).
\end{equation}

%\begin{eqnarray}\label{sys-tele-symdet-1}
%&&\no-\frac{t^{-\alpha}}{\Gamma(1-\alpha)}h_2(x)u v+\big[h_2(x)v+g_1(x)-h_1(x)u-g_2(x)-\xi_x+(\gamma-\alpha) (2\chi_2t+\chi_1)\big]v_x\\
%&& -\bigg[\Big[h_1(x)u_x+\frac{d g_2(x)}{dx}+\frac{d h_1(x)}{dx}u\Big]v+\frac{d f_1(x)}{dx}u+f_1(x)u_x\bigg]+\big[h_2(x)u+f_2(x)\big]\big[P(u)u_x+G(u)\big]=0.
%\end{eqnarray}\widehat{P}

Now we turn to the second equation in system (\ref{tele}). Since system (\ref{tele}) is nonlinear, then either $P(u)$ is a nonconstant function or $Q(u)$ is a nonlinear function or both of them. In order to facilitate symmetry classification, we divide $P(u)$ and $G(u)$ into two parts respectively, set $P(u)=\widehat{P}+\beta$, $G(u)=\lambda\,u+\widehat{G}$, where $\lambda$ and $\beta$ are two constants, $\widehat{P}=\widehat{P}(u)$ is a function without containing constant term and $\widehat{G}=\widehat{G}(u)$ is a nonlinear function without containing the linear term $\lambda\, u$.

With the above assumptions, $I_2=\{u, u_x, \widehat{P}u_x, \widehat{G}\}, J_2=\{u, u_x\}, I_2\setminus J_2=\{\widehat{P}u_x, \widehat{G}\}$ and $\F_2=\widehat{P}u_x+\widehat{G}+\beta u_x+\lambda\, u$. By Theorem \ref{cor-1}, the two determining conditions for the second equation in system (\ref{tele}) are
\begin{subequations}
\begin{equation}\label{sys-tele-3-1}
\frac{\partial^\alpha h_2(t,x)}{\partial t^\alpha}-\beta\frac{\partial  h_1(t,x)}{\partial x}-\lambda\, h_1(t,x)=0,
\end{equation}
\begin{eqnarray}\label{sys-tele-3-2}
&&\no \left(c_1-\alpha\chi_1\right)\big(\widehat{P}u_x+\widehat{G}+\beta u_x+\lambda \,u\big)-\lambda\left[\eta-h_1(t,x)\right]\\
&& \hspace{1.3cm}-\beta\left[\eta^x-\frac{\partial  h_1(t,x)}{\partial x}\right]-\eta\left(\widehat{P}'u_x+\widehat{G}'\right)-\eta^x \widehat{P}=0,
\end{eqnarray}
\end{subequations}
where $\eta^x$ is  {given} by
\begin{equation}
\no \eta^x=\big[ g_1(x)-\xi'(x)\big]u_x+g_1'(x)u+\frac{\partial h_1(t,x)}{\partial x}.
\end{equation}

We claim that if $h_1(t,x)$ is a constant, then $h_1(t,x)=0$. Assume  {that} $h_1(t,x)=C$ is a constant, then solving the first equation in system (\ref{sys-tele-symdet-11}) gives $h_2(t,x)=C xt^{-\alpha}/\Gamma(1-\alpha)+k(t)$. Inserting them into Eq.(\ref{sys-tele-3-1}) yields
\begin{equation}
\no\frac{\partial^\alpha k(t)}{\partial t^\alpha}+C x\frac{\partial^\alpha}{\partial t^\alpha}\left[\frac{t^{-\alpha}}{\Gamma(1-\alpha)}\right]-\lambda\, C=0.
\end{equation}
Equating the  {coefficients} of $x$ to zero yields $C=0$ and then $k(t)=c_2t^{\alpha-1}$. Thus $h_1(t,x)=0,\, h_2(t,x)=c_2t^{\alpha-1}$.

 {Then} inserting $\eta$ and $\eta^x$ into Eq.(\ref{sys-tele-3-2}) and separating it with respect to $u_x$ yields
\begin{subequations}\label{sys63}
\begin{eqnarray}\label{sep-2}
&& 2(\widehat{P}+\beta) \left[c_1-g_1(x)\right]-\widehat{P}' \left[g_1(x)u +h_1(t,x)\right]=0,
\end{eqnarray}
\begin{eqnarray}\label{sep-1}
&&\no \widehat{G} \left(c_1-\alpha \chi_1\right)+u
   \left[\lambda c_1-\lambda g_1(x)-\beta g_1'(x)-\lambda \alpha \chi_1\right] \\
&&\hspace{1.3cm}-\widehat{P} \left[g_1'(x)u+\frac{\partial  h_1(t,x)}{\partial x}\right]-\widehat{G}' \left[g_1(x)u +h_1(t,x)\right]=0,
\end{eqnarray}
\end{subequations}
where system (\ref{sys-tele-symdet-11}) is used. We start with Eq.(\ref{sep-2}) to classify the pairs $(P(u),G(u))$ and then to determine the corresponding Lie symmetries.
First consider $\widehat{P}\neq 0$. Then Eq.(\ref{sep-2}) implies $g_1(x)$ and $h_1(t,x)$ are constants respectively, and thus $h_1(t,x)=0,\, h_2(t,x)=c_2t^{\alpha-1}$ by the claim. Let $g_1(x)=\omega$ {be a constant}. If $\omega=0$, Eq.(\ref{sep-2}) gives $c_1=0$, then we find that $P(u)$ is arbitrary and $G(u)=0$ from system (\ref{sys63}), which is contradict with nonzero $G(u)$. %Then by solving system (\ref{sys-tele-symdet-11}) together with (\ref{gamma-tele-eta-new}), we obtain
%\begin{equation}
%\no\tau=\chi_1 t,\ \xi=\alpha\chi_1,\  \eta=g_1(x)u+h_1(t,x),\ \phi=c_1v+h_2(t,x).
%\end{equation}
While for $\omega\neq0$, solving Eq.(\ref{sep-2}) gives $\widehat{P}=u^{2 (c_1-\omega)/\omega}-\beta$, which implies $\beta=0$ since $\widehat{P}$ is independent of constant term, i.e. $\widehat{P}=u^{2 (c_1-\omega)/\omega}$.

 {Furthermore,} from the second equation in system (\ref{sys-tele-symdet-11}), we find $\xi (x)=(\alpha\chi_1-\omega+c_1)x+c_3$ and then Eq.(\ref{sep-1}) becomes
\begin{eqnarray}\label{sep-6}
&& \widehat{G} \left(c_1-\alpha \chi_1\right)+\lambda u
   \left( c_1-\omega-\alpha \chi_1\right) - \omega u \widehat{G}'=0,
\end{eqnarray}
which gives $ \widehat{G} =u{}^{(c_1-\alpha  \chi_1)/\omega}$ and $\lambda=0$ since $\widehat{G}$ does not contain the term $u$. Therefore, system (\ref{tele}) with $P(u)=u^{2 c_1/\omega-2},G(u) =u{}^{(c_1-\alpha  \chi_1)/\omega}$  {has a Lie symmetry with the infinitesimal generator}
\begin{eqnarray}
&&\no X=\chi_1t\partial_t+\left[(\alpha\chi_1-\omega+c_1)x+c_3\right]\partial_x+\omega u\partial_u+\left(c_1v+c_2t^{\alpha-1}\right)\partial_v.
\end{eqnarray}

If $\widehat{P}=0$, then $\beta\neq0$ since $P(u)\neq0$. We find $g_1(x)=c_1$ from Eq.(\ref{sep-2}) and $h_1(t,x)=0$ by the claim and Eq.(\ref{sep-1}). Under such conditions, Eq.(\ref{sep-1}) becomes Eq.(\ref{sep-6}) with $\omega=c_1$, which is a particular case of $\widehat{P}\neq 0$.

We summarize the above Lie symmetry classifications of system (\ref{tele}) as the following proposition.
\begin{prop}\label{ro0-1}
The Lie symmetries admitted by the time-fractional nonlinear telegraph equations (\ref{tele}) are classified as follows:

I). For arbitrary functions $P(u)$ and $G(u)$, system (\ref{tele}) is admitted by $X=\partial_x$.

II). For $P(u)=u^{2 c_1/\omega-2}, G(u) =u^{(c_1-\alpha  \chi_1)/\omega}$, system (\ref{tele}) is admitted by
\begin{eqnarray}
&&\no X=\chi_1t\partial_t+\left[(\alpha\chi_1-\omega+c_1)x+c_3\right]\partial_x+\omega u\partial_u+\left(c_1v+c_2t^{\alpha-1}\right)\partial_v.
\end{eqnarray}
\end{prop}

Note that the cases of symmetry classification for fractional {nonlinear telegraph equations} (\ref{tele}) decrease significantly compared with the ones of integer-order case in \cite{chao}, the reason for such phenomenon is that the fractional derivative greatly affects the symmetry properties of fractional PDEs.

With the infinitesimal generators in Proposition \ref{ro0-1}, we perform symmetry reductions for system (\ref{tele}) as the following three cases.

I). For arbitrary functions $P(u)$ and $G(u)$, system (\ref{tele}) is admitted by the infinitesimal generator $X=\partial_x$. Similar as the case for system (\ref{kdv}), we obtain a solution of system (\ref{tele}) $u(t,x)=C_1t^{\alpha-1},v(t,x)=C_2t^{\alpha-1}$, where $C_1$ { and} $C_2$ are integral constants.

II). For $P(u)=u^{2 c_1/\omega-2}, G(u) =u^{(c_1-\alpha \chi_1)/\omega}$, induced by the infinitesimal generator $X=\partial_x+ c_2t^{\alpha-1}\partial_v$, system (\ref{tele}) is {reduced to}
\begin{eqnarray}\label{reud}
&& \frac{\partial^\alpha U(\zeta)}{\partial\zeta^\alpha}=c_2\zeta^{\alpha-1},~~~~\frac{\partial^\alpha V(\zeta)}{\partial \zeta^\alpha}=U(\zeta)^{\frac{1}{\omega}(c_1-\alpha \chi_1)},
\end{eqnarray}
where the similarity variables are $\zeta=t, \, U(\zeta)=u(t,x)$ { and} $V(\zeta)=v(t,x)-x t^{\alpha-1}$. Then {by solving the reduced equations (\ref{reud}) and changing to original variables}, we obtain a particular solution of system (\ref{tele})
\begin{eqnarray}
&&\no u(t,x)=\frac{c_2\Gamma(\alpha)}{\Gamma(2\alpha)}t^{2\alpha-1},~~~ v(t,x)=x t^{\alpha-1}+ \left[\frac{c_2\Gamma(\alpha)}{\Gamma(2\alpha)}\right]^{\vartheta} \frac{\Gamma((2\alpha-1)\vartheta+1)}{\Gamma((2\alpha-1)\vartheta+\alpha+1)}t^{(2\alpha-1)\vartheta+\alpha},
\end{eqnarray}
where $\vartheta=(c_1-\alpha \chi_1)/\omega$.

III). For $P(u)=u^{2 c_1/\omega-2}, G(u) =u^{(c_1-\alpha \chi_1)/\omega}$, by the infinitesimal generator $X=\chi_1 t \partial_t+(c_1+\chi_1\alpha-\omega) x \partial_x+ \omega u\partial_u+ c_1v\partial_v$, system (\ref{tele}) is reduced {to}
\begin{eqnarray}
&&\no \left(P_{\frac{\chi_1}{c_1+\chi_1\alpha-\omega}}^{1+\frac{\omega}{\chi_1}-\alpha,\,\alpha}U\right)(\zeta)=V'(\zeta),\\
&&\no\left(P_{\frac{\chi_1}{c_1+\chi_1\alpha-\omega}}^{1+\frac{c_1}{\chi_1}-\alpha,\,\alpha}V\right)(\zeta)=U(\zeta)U'(\zeta)+U^{(c_1-\alpha\chi_1)/\omega}(\zeta),
\end{eqnarray}
where $\zeta=xt^{-(c_1+\alpha \chi_1-\omega)/\chi_1}, \, U(\zeta)=u(t,x)t^{-\omega/\chi_1}$ { and} $ V(\zeta)=v(t,x)t^{-c_1/\chi_1}$.

\section{Conclusions}
{We deeply study the Lie group theory for the system of multi-dimensional time-fractional PDEs and give an explicit formula of the extended infinitesimal involving Riemann-Liouville fractional derivative. Moreover, we show that {the infinitesimal generators of Lie symmetries for system  (\ref{eqn-ger}) have an elegant structure and are} completely determined by two elegant conditions. Our results pave a simple way for acquiring symmetry information of multi-dimensional time-fractional PDEs (\ref{eqn-ger})  and also motivate the automatic implementation of searching for Lie symmetries with lower complexity of computation.}% Applications to three examples demonstrate the efficiency of the results. Our results revise the mistakes in the previous literatures though they do not affect the obtainment of Lie symmetries of time-fractional PDEs.

%It is worthy of saying thatIn addition, whether the above results hold for time-fractional PDEs with Caputo fractional derivative are still unknown.  We will proceed in this way and report the related results.pave a way for that the whole procedure is mechanic and  performed on the computer , and thus
\bigskip

\section*{Acknowledgements}
This paper is supported by the National Natural Science Foundation of China (No. 11671014).
\\\\
\textbf{Declarations of interest: none}

\end{document}